\documentclass[a4paper,11pt]{article}
\usepackage{amsmath}
\usepackage{amssymb}
\usepackage[dvipdfmx]{graphicx,color}
\nopagebreak[3]
\newcommand{\newsection}[1]
{\section{#1}\setcounter{theorem}{0} \setcounter{equation}{0} \par\noindent}

\newtheorem{theorem}{Theorem}

\newtheorem{lemma}[theorem]{Lemma}

\newtheorem{proposition}[theorem]{Proposition}

\newtheorem{claim}[theorem]{Claim}
\newcommand{\beq}{ \begin{equation} }
\newcommand{\eeq}{ \end{equation} }
\newcommand{\ebox}{\hfill$\square$ }

\newcommand{\br}{{\mathbb R}}

\newcommand{\supp}{\mbox{\rm supp}\ }
\renewcommand{\div}{\mbox{\rm div}\ }
\newcommand{\proof}{\noindent{\bf Proof.\ } }
\newcommand{\vep}{ \varepsilon }

\newcommand{\till}{ {\tilde L} }
\newcommand{\anglex}{ {\langle x\rangle} }
\newcommand{\angley}{ {\langle y\rangle} }

\newcommand{\angles}{ {\langle s\rangle} }
\newcommand{\angler}{ {\langle r\rangle} }

\newcommand{\angletpr}{ {\langle t+r\rangle} }
\newcommand{\angletmr}{ {\langle t-r\rangle} }
\newcommand{\anglecItmr}{ {\langle c_It-r\rangle} }

\newcommand{\anglectmr}{ {\langle ct-r\rangle} }

\newcommand{\anglectmy}{ {\langle ct-|y|\rangle} }
\newcommand{\partialbar}{ \overline{\partial} }
\newcommand{\ext}{ { \mathbb{R}^3\backslash \mathcal{K} } }
\newcommand{\whole}{ { \mathbb{R}^3 } }
\newcommand{\Rn}{ { \mathbb{R}^n } }
\title{
Remarks on global solutions for nonlinear wave equations under the standard null conditions
}
\author
{
Hans Lindblad
\thanks
{
Department of Mathematics, 
Johns Hopkins University, 
404 Krieger Hall, 
3400 N. Charles Street,
Baltimore, MD 21218.
\texttt{lindblad@math.jhu.edu}
}
\and
Makoto Nakamura
\thanks
{Mathematical Institute, 
Tohoku University, Sendai 980-8578.
\texttt{makoto@math.tohoku.ac.jp}
}
\and
Christopher D. Sogge
\thanks
{
Department of Mathematics, 
Johns Hopkins University, 
404 Krieger Hall, 
3400 N. Charles Street,
Baltimore, MD 21218.
\texttt{sogge@math.jhu.edu}
}
}

\date{}

\begin{document}

\maketitle

\begin{abstract}
A combination of some weighted energy estimates is applied for 
the Cauchy problem of quasilinear wave equations with the standard null conditions in three spatial dimensions.
Alternative proofs for global solutions are shown including the exterior domain problems.
\end{abstract}
{\small 
\begin{tabular}{l}
\end{tabular}
} 
%
\section{Introduction}

Let us consider the Cauchy problem of nonlinear wave equation under the standard null conditions 
\begin{equation}
\label{Cauchy-Intro}
(\partial_t^2-\Delta )u(t,x)=(\partial_t u(t,x))^2-|\nabla u(t,x)|^2
\ \ \ \ \mbox{for}\ \ (t,x)\in [0,\infty)\times\whole,
\end{equation}
where $\Delta:=\partial_1^2+\partial_2^2+\partial_3^2$, and $\nabla:=(\partial_1,\partial_2,\partial_3)$.
In \cite[p52]{Lindblad-Rodnianski-2005-CMP},  it is pointed out that the combination of
the weighted energy estimate 
\begin{multline}
\label{Intro-LR}
\|\nabla_{t,x}u\|_{L^\infty((0,T),L^2(\mathbb{R}^3))}^2
+
\|\langle t-r\rangle^{-1/2-\gamma} \partialbar u\|^2
_{L^2((0,T)\times\mathbb{R}^3)} \\
\le
C\|\nabla_{t,x}u(0,\cdot)\|_{L^2(\mathbb{R}^3)}^2
+C\int_0^T \int_\whole
|(\partial_tu)(\partial_t^2-\Delta)u|
dxdt
\end{multline}
(see \cite[p76, Corollary 8.2]{Lindblad-Rodnianski-2005-CMP} 
and Alinhac \cite[Theorem 1]{Alinhac-2004-MathAnn}), 
where $T>0$, $r=|x|$, $\nabla_{t,x}=(\partial_t, \partial_1,\partial_2,\partial_3)$,  
$\partialbar=(\partial_t+\partial_r, \nabla-\frac{x}{r}\partial_r)$ 
denotes tangential derivatives along the light cone $t=r$, the constant $C$ is independent of $T$, 
$\gamma>0$ is any fixed real number, 
and the Klainerman-Sobolev estimate 
\begin{equation}
(1+t+r)(1+|t-r|)^{1/2}|u(t,x)|
\le 
C \sum_{|\alpha|\le 2} \|\Theta^\alpha u(t,\cdot)\|_{L^2(\mathbb{R}^3)}
\end{equation}
(see \cite[p118, Proposition 6.5.1]{Hoermander-1997-Springer}) 
gives a much simplified proof for the existence of the small solutions for \eqref{Cauchy-Intro}, 
where $\Theta$ denotes the vector fields 
\begin{equation}
\partial_t,\ \ \partial_j,\ \ 
t\partial_t+rdr,\ \ 
t\partial_j+x_j\partial_t,\ \ 
x_j\partial_k-x_k\partial_j, 
\ \ 1\le j\neq k\le 3,
\end{equation}
and $\alpha$ denotes multiple indices.
In this paper, we generalize this argument so that we are able to treat the system of wave equations with different speeds and also the corresponding exterior domain problems.
While the above weighted energy estimates could be generalized to treat the $c>0$ speed D'Alembertian $\square_c:=\partial_t^2-c^2\Delta$, the operators $\{t\partial_j+x_j\partial_t\}_{j=1}^3$ are not commutable with $\square_c$ if $c\neq1$, and they also makes it difficult to handle the energy near the obstacle when we consider the exterior domain problems.
To avoid the use of these operators, we use the Klainerman-Sobolev type estimates by Sideris, Tu, Hidano and Yokoyama 
\begin{multline}
\label{Intro-KS}
\langle r\rangle^{1/2} \langle r+t \rangle^{1/2}
\langle ct-r\rangle^{1/2}
|u'(t,x)|
\le 
C\sum_{\substack{ \mu+|\alpha|\le2 \\ \mu \le 1} } 
\|L^\mu Z^\alpha u'(t,\cdot)\|_{L^2(\mathbb{R}^3)}
\\
+
C\sum_{ \substack{|\alpha|\le1 } } 
\|\langle t+|y|\rangle Z^\alpha (\partial_t^2-c^2\Delta) u(t,y)\|_{L^2_y(\mathbb{R}^3)},
\end{multline}
where $L:=t\partial_t+r\partial_r$, 
$Z:=(\partial_t, \nabla, \{x_j\partial_k-x_k\partial_j\}_{1\le j\neq k\le3})$
(see Lemma \ref{Lemma-KS}, below). 
And we also use the weighted energy estimate by Keel, Smith and Sogge 
\begin{multline}
\label{Intro-KSS}
\|\nabla_{t,x}u\|_{L^\infty((0,T),L^2(\mathbb{R}^3))}
+
(\log(e+T))^{-1/2}\|\anglex^{-1/2}\nabla_{t,x}u(t,x)\|
_{L^2((0,T)\times\mathbb{R}^3)}
\\
\le
C\|\nabla_{t,x} u(0,\cdot)\|_{L^2(\mathbb{R}^3)}
+
C\|(\partial_t^2-c^2\Delta)u\|_{L^1((0,T),L^2(\mathbb{R}^3))}
\end{multline}
(see Lemma \ref{Lemma-EnergyKSSTangential}, below).
We remark that the combination of \eqref{Intro-LR}, \eqref{Intro-KS} and \eqref{Intro-KSS} gives a much  simplified and elementary proof for the global existence for the small solutions of nonlinear wave equations under the standard multispeed null conditions. 
Especially, we give alternative proofs of the following two theorems.

%

\subsection{The Cauchy problems without obstacles}

We consider the Cauchy problem for a system of quasilinear wave equations with $D\ge1$ propagation speeds $\{c_I\}_{1\le I\le D}$, $c_I>0$.
We put $u=(u_1,\cdots,u_D)$, $F=(F_1,\cdots,F_D)$, $f=(f_1,\cdots,f_D)$, $g=(g_1,\cdots,g_D)$, and we consider  
\begin{equation}
\label{P-Whole}
\left\{
\begin{array}{l}
(\partial_t^2-c_I^2\Delta
)u_I(t,x)=F_I
(u', u'')(t,x)
\ \ \mbox{for}\ t\in [0,\infty),\ x\in \whole,\ 1\le I\le D \\
u(0,\cdot)=f(\cdot),\ \ \partial_tu(0,\cdot)=g(\cdot),
\end{array}
\right.
\end{equation}
where we put $\partial_0=\partial_t$, we denote 
the first derivatives $\{\partial_j u\}_{0\le j\le 3}$ by $u'$, and the second derivatives $\{\partial_j\partial_k u\}_{0\le j,k\le 3}$ by $u''$.
We assume that $F$ vanishes to the second order and has the form 
\begin{equation}
\label{F}
F_I(u',u'')=B_I(u')+Q_I(u', u''),
\end{equation}
where 
\begin{equation}
\label{B}
B_I(u'):= 
\sum_{\substack{1\le J,K\le D \\ 0\le j,k\le 3} }
B_I^{JKjk}\partial_j u_J\partial_k u_K
\end{equation}
\begin{equation}
\label{Q}
Q_I(u',u''):= 
\sum_{\substack{1\le J,K\le D \\ 0\le j,k,l\le 3} }
Q_I^{JKjkl}\partial_j u_J\partial_k \partial_l u_K.
\end{equation}
Here, 
$\{B_I^{JKjk}\}_{\substack{1\le I,J,K\le D \\ 0\le j,k\le 3} }$ and 
$\{Q_I^{JKjkl}\}_{\substack{1\le I,J,K\le D \\ 0\le j,k,l\le 3} }$ are real constants
which satisfy the symmetry condition
\begin{equation}
\label{Symmetry}
Q_I^{JKjkl}=Q_K^{JIjlk},
\end{equation}
which is used to derive the energy estimates.
To show the global solutions, we assume the standard null conditions 
\begin{equation}
\label{Null-Condition}
\sum_{0\le j,k\le 3} B_I^{JKjk} \xi_j\xi_k=\sum_{0\le j,k,l\le 3} Q_I^{JKjkl} \xi_j\xi_k\xi_l=0
\end{equation}
for any $1\le I, J, K\le D$ with $c_I=c_J=c_K$, and any 
$(\xi_0,\xi_1,\xi_2, \xi_3)\in \mathbb{R}^4$ with   
$\xi_0^2=c_I^2(\xi_1^2+\xi_2^2+\xi_3^2)$.
For  example, 
\beq
B_I(u')=
\sum_{\substack{1\le J,K\le D \\ (c_J,c_K)=(c_I,c_I)} }
\kappa^{JK} \{(\partial_t u_J)(\partial_t u_K)-c_I^2 (\nabla u_J)\cdot(\nabla u_K) \}
+
\sum_{\substack{1\le J,K\le D \\ (c_J,c_K)\neq(c_I,c_I)} }
\lambda^{JK} u_J'u_K'
\eeq
for $\{ \kappa^{JK} \}_{1\le J,K\le D}$, $\{\lambda^{JK}\}_{1\le J,K\le D}\subset \mathbb{R}$,  and $Q_I(u',u'')=\partial_{t,x} B_I(u')$ satisfy the null conditions.
It is known that any nontrivial solutions blow up in finite time in general for quadratic nonlinearities
(see John \cite{John-1981-CPAM}), 
while the null conditions guarantee the global solutions
(see Christodoulou 
\cite{Christodoulou-1986-CPAM}
 and Klainerman 
\cite{Klainerman-1986-LecturesApplMath}).
We give an alternative proof of the following theorem, which has been shown by Sideris and Tu in 
\cite{Sideris-Tu-2001-SIAM}.

\begin{theorem}
\label{Theorem}
Let $f$ and $g$ be smooth functions, and let $F$ satisfy the above standard null conditions.
Then there exists a positive natural number $N$ 
(for example, we are able to take $N=12$)
such that if 
\begin{equation}
\sum_{|\alpha|\le N}\|\anglex^{|\alpha|}\partial_x^\alpha \nabla f(x)\|_{L^2(\whole)}
+
\sum_{|\alpha|\le N}\|\anglex^{|\alpha|}\partial_x^\alpha g(x)\|_{L^2(\whole)}
\end{equation}
is sufficiently small, then there exists a unique global solution $u\in C^\infty([0,\infty)\times\whole)$ of \eqref{P-Whole}.
\end{theorem}

The proof in \cite{Sideris-Tu-2001-SIAM} consists of the standard energy estimates and a series of pointwise estimates for 
$r^{1/2}u$, 
$\langle r\rangle u'$, 
$\langle r\rangle \langle ct-r \rangle^{1/2} u'$,  
$\langle r\rangle \langle ct-r \rangle u''$, 
and the weighted estimate 
$\|\langle ct-r \rangle u''\|_{L^2(\mathbb{R}^3)}$ 
plays an important role to control the energy near the light cone $r=ct$.
The weighted energy estimates in \eqref{Intro-LR} and \eqref{Intro-KSS} are not used in \cite{Sideris-Tu-2001-SIAM}.
We remark in this paper the combination of \eqref{Intro-LR} and \eqref{Intro-KSS} yields much simplified and elementary proof for the theorem.
The key estimate is the lower energy estimate \eqref{Ineq-LowEnergy}, 
which is proved via a straightforward application of
\eqref{Intro-LR}, \eqref{Intro-KS}, \eqref{Intro-KSS} and the estimate for null conditions Lemma \ref{Lemma-NullConditions}.
We generalize \eqref{Intro-LR} in Lemma \ref{Lemma-WE}.
It is interesting to see that Lemma \ref{Lemma-WE} has a close similarity to \eqref{Intro-KSS} 
(cf. \cite[Corollary 5]{Nakamura-2011-DIE}).
The estimate \eqref{Intro-KSS} yields much simplified proof for almost global solutions 
(\cite{Keel-Smith-Sogge-2002-JAnalMath}) and also has strong applications to global solutions 
(\cite{Metcalfe-Nakamura-2007-JDE, Metcalfe-Nakamura-Sogge-2005-ForumMath, Metcalfe-Nakamura-Sogge-2005-JJM, Metcalfe-Sogge-2005-InventMath, Nakamura-2012-JDE}).

%
\subsection{Exterior domain problems}

We also consider the exterior domain problem. 
Let $\mathcal{K}$ be any fixed compact domain in 
$\mathbb{R}^3$ with smooth boundary.
Without loss of generality, we may assume $0\in \mathcal{K}\subset \{x\in \mathbb{R}^3 : |x|<1\}$ by the shift and scaling arguments.
Moreover, we assume the following local energy decay estimates.
Let $u$ be the solution of 
\begin{equation}
\left\{
\begin{array}{l}
(\partial_t^2-\Delta) u(t,x)=0\ \ \ \ \mbox{for}\ \ (t,x)\in [0,\infty)\times \mathbb{R}^3\backslash \mathcal{K} \\
u(t,\cdot)|_{\mathcal{K}} =0\ \ \ \ \mbox{for}\ \ t\in [0,\infty) \\
u(0, \cdot)=f(\cdot),\ \ \partial_tu(0,\cdot)=g(\cdot).
\end{array}
\right.
\end{equation}
If the initial data satisfy $\supp f\cup\supp g\subset \{x\in \ext : |x|<4\}$, then there exist constants $C>0$ and $a>0$ such that 
\begin{equation}
\label{LED}
\|u'(t,\cdot)\|_{ L^2(\{x\in \ext : |x|<4\}) }
\le Ce^{-at} \sum_{|\alpha|\le1} \|\partial_x^\alpha u'(0,\cdot)\|_{L^2(\ext)}
\end{equation}
for any $t\ge0$.
The local energy decay estimate \eqref{LED} holds if the obstacle is nontrapping without the loss of derivatives $|\alpha|=1$ 
(Morawetz, Ralston and Strauss 
\cite{Morawetz-Ralston-Strauss-1977-CPAM}),
or the obstacle consists of certain finite unions of convex obstacles 
(Ikawa 
\cite{Ikawa-1982-OsakaJMath, Ikawa-1998-AnnInstFourier}).

We consider the exterior domain problems for \eqref{P-Whole} given by 
\begin{equation}
\label{P}
\left\{
\begin{array}{l}
(\partial_t^2-c_I^2\Delta)u_I(t,x)=F_I(u', u'')(t,x)
\ \ \mbox{for}\ t\in [0,\infty),\  x\in \ext,\ 1\le I\le D \\
u(t,x)|_{x\in \partial\mathcal{K} } =0
\ \ \mbox{for}\ \ t\in [0,\infty) \\
u(0,\cdot)=f(\cdot),\ \ \partial_tu(0,\cdot)=g(\cdot),
\end{array}
\right.
\end{equation}
where $F_I$ is written as \eqref{F}, \eqref{B}, \eqref{Q}, we assume the symmetry condition  \eqref{Symmetry} and the null conditions \eqref{Null-Condition}.
Since \eqref{P} is the initial and boundary value problem, the initial data $f$ and $g$ must satisfy the compatibility condition.
For $k\ge0$ and the solution $u$ of \eqref{P},  the condition $\partial_t^k u(0,\cdot)=0$ is written in terms of $f$, $g$ and $F$.
We assume the compatibility condition of infinite order, namely, 
$\partial_t^k u(0,\cdot)|_{\mathcal{K} }=0$ for any $k\ge0$.
We give an alternative proof of the following theorem, 
which has been shown in 
\cite{Metcalfe-Nakamura-Sogge-2005-ForumMath, Metcalfe-Sogge-2005-InventMath}.

%
\begin{theorem}
\label{Theorem-Ext}
Let $f$ and $g$ be smooth functions and satisfy the compatibility conditions of infinite order.
Let $F$ satisfy the above standard null conditions.
Then there exists a positive large natural number $N$ 
(for example, we are able to take $N=64$) 
such that if 
\begin{equation}
\sum_{|\alpha|\le N}\|\anglex^{|\alpha|+1}\partial_x^\alpha f(x)\|_{L^2(\ext)}
+
\sum_{|\alpha|\le N-1}\|\anglex^{|\alpha|+1}\partial_x^\alpha g(x)\|_{L^2(\ext)}
\end{equation}
is sufficiently small, then \eqref{P} has a unique global solution 
$u\in C^\infty([0,\infty)\times\ext)$.
\end{theorem}

There is a series of papers on almost global and global solutions 
by Keel, Smith and Sogge
\cite{Keel-Smith-Sogge-2000-AmerJMath} for convex obstacles,
\cite{Keel-Smith-Sogge-2002-JAnalMath} for nontrapping obstacles,
\cite{Keel-Smith-Sogge-2002-JFA}  
and \cite{Keel-Smith-Sogge-2004-JAMS} for star-shaped obstacles.
See also 
\cite{Metcalfe-Nakamura-Sogge-2005-ForumMath},
\cite{Metcalfe-Nakamura-Sogge-2005-JJM} 
and 
\cite{Metcalfe-Sogge-2005-InventMath} for Ikawa's type trapping obstacles.


In \cite{Metcalfe-Nakamura-Sogge-2005-ForumMath, Metcalfe-Sogge-2005-InventMath}, 
the weighted estimate \eqref{Intro-KSS} has been used, while tangential derivatives and \eqref{Intro-LR} are not used.
Let $B(u,v)=\sum_{0\le j,k\le 3}B^{jk}\partial_ju\partial_kv$ 
satisfy the null condition 
$\sum_{0\le j,k\le 3}B^{jk}\xi_j\xi_k=0$ 
for $\xi_0^2=\sum_{j=1}^3\xi_j^3$.
One of the advantages to use the tangential derivatives is that we are able to estimate the null conditions simply as 
\beq
\label{Intro-Estimates-Null}
| \sum_{\substack{0\le j,k\le 3} }
B^{jk} 
\partial_j u\partial_k v |
\le 
C|\partialbar u| |v'|
+
C|u'| |\partialbar v|
\eeq
(see Lemma \ref{Lemma-NullConditions}). 
In \cite{Metcalfe-Nakamura-Sogge-2005-ForumMath,
Metcalfe-Sogge-2005-InventMath, 
Sideris-Tu-2001-SIAM},  
the type of estimate 
\begin{multline}
|\sum_{0\le j,k\le 3}B^{jk} \partial_j u\partial_k v|
\le 
\frac{C}{\langle r\rangle}
\cdot
\left\{
\sum_{\substack{\mu+|\alpha|\le 1 \\ \mu\le 1} }
|L^\mu Z^\alpha u||\partial v|+|\partial u|
\sum_{\substack{\mu+|\alpha|\le 1 \\ \mu\le 1} }
|L^\mu Z^\alpha v|
\right\}
\\
+
C\frac{\langle t-r\rangle}{\langle t+r\rangle}
|\partial u||\partial v|
\end{multline}
(see \cite[Lemma 5.4]{Sogge-2008-InterPress}) 
has been used, which needs variants of Sobolev type estimates, or $L^\infty-L^1$ estimates based on 
the Kirchhoff formula 
to bound $|L^\mu Z^\alpha u|$ and $|L^\mu Z^\alpha v|$, 
and the proof for global solutions needs structural complexity.
In this paper, we use \eqref{Intro-Estimates-Null}, and we show the combination of two type of weighted energy estimates \eqref{Intro-LR}, \eqref{Intro-KSS}, and the Sobolev estimates \eqref{Intro-KS} gives a simplified proof of the theorem.

\subsection{Notation}
We use the method of commuting vector fields introduced by John and Klainerman 
\cite{John-1983-CPAM, John-Klainerman-1984-CPAM, Klainerman-1984-CPAM}.
See also Keel, Smith and Sogge 
\cite{Keel-Smith-Sogge-2002-JAnalMath} for exterior domains.
We denote the space-time derivatives by $\partial$, the rotational derivatives by $\Omega$, and the scaling operator by $L$.
We use $u'$ to denote $\partial u$ in some cases.
We denote $\partial, \Omega$ by $Z$, and $\partial, \Omega, L$ by $\Gamma$.
We use the tangential derivatives $\partialbar_c$ on the $c$-speed light cone $ct=|x|$ to treat the nonlinear terms which satisfies the null conditions.
We summarize as 
\begin{equation}
\left\{
\begin{array}{l}
x=r\omega, 
\ \ 
\omega\in \mathbb{S}^{2}, 
\ \ 
x_0=t,
\ \ \partial_0=\partial_t,
\\ 
\partial=(\partial_0,\partial_1,\partial_2,\partial_3),
\ \ 
\Omega=(x_j\partial_k-x_k\partial_j)_{1\le j\neq k\le 3},
\ \ 
L=t\partial_t+r\partial_r,
\\ 
Z=(\partial, \Omega),\ \ \Gamma=(\partial, \Omega, L),
\\
\partialbar_c=(\partialbar_{c0},\partialbar_{c1},\partialbar_{c2},\partialbar_{c3})
=(\partial_t+c\partial_r, \nabla_x-\omega\partial_r) \mbox{ for }c>0.
\end{array}
\right.
\end{equation}
The operators $Z$, $L$ have the commuting properties with the  $c$ speed D'Alembertian $\square_c:=\partial_t^2-c^2\Delta$ such as 
\begin{equation}
\square_cZ=Z\square_c,\ \ \ \ \square_cL=(L+2)\square_c.
\end{equation}
We do not use the Lorentz boosts $\{t\partial_j+x_j\partial_t\}_{j=1}^3$ which are not suitable for the different speeds system or the exterior domain problems.
We put $\angler=\sqrt{1+r^2}$, $\anglex=\sqrt{1+|x|^2}$.
The Lebesgue spaces on $\whole$ are denoted by $L^p(\whole)$ for $1\le p\le \infty$. 
For $R>0$, $L^p(|x|<R)$ denotes $L^p(\{x\in \whole : |x|<R\})$.
The strip region is denoted by $S_T:=[0,T)\times \whole$ for $T\ge 0$.
We use $L^p(\ext)$, $L^p(|x|<R):=L^p(\{x\in \ext : |x|<R\})$ and $S_T:=[0,T)\times (\ext)$ 
when we consider the exterior domain problems.
For $c>0$, $\square_c:=\partial_t^2-c^2\Delta$ denotes the $c$-speed D'Alembertian.
We put $\square:=\square_1=(\partial_t^2-\Delta)$.
Throughout this paper, $C$ denotes a positive constant which may differ from line to line.
The notation $a\lesssim b$ denotes the inequality $a\le Cb$ for a positive constant $C$ which is not essential for our arguments.
This paper is organized as follows.
In sections 2 and 4, we prepare several estimates which are needed to show Theorem \ref{Theorem} and \ref{Theorem-Ext}, respectively.
Theorems \ref{Theorem} and \ref{Theorem-Ext} are shown in sections 3 and 5, respectively.
In section 6, we put two appendices.
The first is for the proof of the weighted energy estimates, and the second is for a remark on the two spatial dimensions.

%

\newsection{Several estimates to prove Theorem \ref{Theorem} }
In this section, we prepare several estimates to prove Theorem \ref{Theorem}.

\subsection{Energy estimates}
Let us consider the general dimension $n\ge1$ in this subsection.
We put $\Delta=\sum_{j=1}^n\partial_j^2$.
We use the following energy estimates for quasilinear wave equations.
Let $\gamma_{I}^{Kkl}$, $1\le I,K\le D$, $0\le k,l\le n$, be functions which satisfy the symmetry conditions 
$
\gamma_{I}^{Kkl}=\gamma_{K}^{Ilk}.
$
We put 
\begin{equation}
\label{Eqn-Energy-Quasi}
\square_{\gamma_I} u_I
:=
(\partial_t^2-c_I^2\Delta)u_I
-
\sum_{\substack{1\le K\le D\\ 0\le k,l\le n}}
\gamma_I^{Kkl} \partial_k\partial_l u_K \ \ \ \ \mbox{for}\ \ 1\le I\le D.
\end{equation}
We define the energy momentums $e_k(u)$, $0\le k \le n$, and the remainder term $R(u)$ by 
\begin{eqnarray*}
e_0(u)&:=&\sum_{1\le I\le D} \{(\partial_tu_I)^2+c_I^2|\nabla u_I|^2\}
-
\sum_{\substack{1\le I,K\le D \\ 0\le l\le n} }
2\gamma_I^{K0l}\partial_0u_I\partial_l u_K 
\\
&& 
+
\sum_{\substack{1\le I,K\le D \\ 0\le k,l\le n} }
\gamma_{I}^{Kkl}\partial_k u_I\partial_l u_K 
\\
e_k(u)&:=&-\sum_{1\le I\le D} 2c_I^2\partial_0u_I\partial_ku_I
-
\sum_{\substack{1\le I,K\le D \\ 0\le l\le n} }
2\gamma_{I}^{Kkl} \partial_0u_I\partial_l u_K
\ \ \ \ \mbox{for}\ \ 1\le k\le n
\\
R(u)&:=&-\sum_{\substack{1\le I,K\le D \\ 0\le k,l\le n} }
2(\partial_k \gamma_{I}^{Kkl})\partial_0u_I\partial_l u_K
+
\sum_{\substack{1\le I,K\le D \\ 0\le k,l \le 3} }
(\partial_0\gamma_{I}^{Kkl})\partial_k u_I\partial_l u_K.
\end{eqnarray*}
Then the multiplication of $2\partial_tu_I$ to the equation \eqref{Eqn-Energy-Quasi} yields the divergence form
\begin{equation}
\label{Eqn-Divergence}
\partial_t e_0(u)+\div(e_1(u),\cdots,e_n(u))
=\sum_{1\le I\le D} 2\partial_t u_I\square_{\gamma_I} u_I+R(u).
\end{equation}

%

\subsection{Weighted energy estimates}

We use the following weighted energy estimates.

\begin{lemma}
\label{Lemma-EnergyKSSTangential}
Let $n\ge1$.
We put $\Delta:=\sum_{j=1}^n\partial_j^2$.
Let $c>0$ and $T>0$.
The solution $u$ of the Cauchy problem
\begin{equation}
\label{Cauchy-Linear}
\left\{
\begin{array}{l}
(\partial_t^2-c^2\Delta) u (t,x)= F(t,x) \ \ \ \ \mbox{for}\ \ 
(t,x)\in [0,T)\times \mathbb{R}^n \\
u(0,\cdot)=f(\cdot), \ \ 
\partial_t u(0,\cdot) =g(\cdot).
\end{array}
\right.
\end{equation}
satisfies the following estimate, where $w'=\partial_{t,x}w$, 
$\partialbar_c=(\partial_t+c\partial_r, \nabla-\frac{x}{r}\partial_r)$.
\begin{multline}
\label{Ineq-EnergyKSSTangential}
\|u'\|_{L^\infty((0,T),L^2(\mathbb{R}^n))}
+
(\log(e+T))^{-1/2}
\|\anglex^{-1/2} u'(\cdot,x)\|_{L^2((0,T)\times\mathbb{R}^n)}
\\
+
(\log(e+T))^{-1/2}
\|\langle ct-r\rangle^{-1/2} \partialbar_c u(t,x)\|
_{L^2((0,T)\times\mathbb{R}^n)} \\
\lesssim
\|\nabla f\|_{L^2(\mathbb{R}^n)}
+
\|g\|_{L^2(\mathbb{R}^n)}
+
\|F\|_{L^1((0,T),L^2(\mathbb{R}^n))}
\end{multline}
\end{lemma}

The estimate for the first term in the left hand side is the standard energy estimate.
The estimate for the second term is due to Keel, Smith and Sogge 
\cite[Proposition 2.1]{Keel-Smith-Sogge-2002-JAnalMath} for $n=3$,  
Metcalfe, Sogge and Hidano 
for general dimensions
(see \cite[Corollary 5]{Nakamura-2011-DIE}).
The estimate for the third term 
is a logarithmic version of the estimate due to Lindblad and Rodnianski 
\cite[p76, Corollary 8.2]{Lindblad-Rodnianski-2005-CMP} who treat the case $n=3$ and $c=1$. 
We generalize it as Lemma \ref{Lemma-WE} in Appendices. 

%

\subsection{Klainerman-Sobolev type estimates}

We also use the following Klainerman-Sobolev type estimates.

\begin{lemma}
\label{Lemma-KS}
Let $c>0$.
The following inequality holds for any function $u$. 
\begin{multline}
\langle r\rangle^{1/2} \langle r+t \rangle^{1/2}
\langle ct-r\rangle^{1/2}
|u'(t,x)|
\lesssim
\sum_{\substack{ \mu+|\alpha|\le2 \\ \mu \le 1} } 
\|L^\mu Z^\alpha u'(t,\cdot)\|_{L^2(\mathbb{R}^3)}
\\
+
\sum_{ \substack{|\alpha|\le1 } } 
\|\langle t+|y|\rangle Z^\alpha (\partial_t^2-c^2\Delta) u(t,y)\|_{L^2_y(\mathbb{R}^3)}.
\end{multline}
\end{lemma}

\proof 
This lemma directly follows from the combination of 
\begin{multline}
\angler \anglectmr^{1/2} |u(t,x)|
\lesssim
\sum_{|\alpha|\le 2} \|\Omega^\alpha u(t,y)\|_{L^2(|y|> r)}
\\
+
\sum_{|\alpha|\le 1} \|\anglectmy\partial_r \Omega^\alpha u(t,y)\|_{L^2(|y|> r)}
\end{multline}
by Sideris \cite[Lemma 3.3]{Sideris-2000-AnnalsMath}, 
\begin{multline}
\angler^{1/2} \anglectmr |u'(t,x)|
\lesssim
\sum_{|\alpha|\le 1} \|Z^\alpha u'(t,\cdot)\|_{L^2(\mathbb{R}^3)}
\\
+
\sum_{\substack{|\alpha|\le 1, |\beta|=2 } } 
\|\anglectmy Z^\alpha \partial^\beta u(t,y)\|_{L^2_y(\mathbb{R}^3)}
\end{multline}
by Hidano \cite[Lemma 4.1]{Hidano-2004-TohokuMathJ}, and 
\begin{multline}
\sum_{ \substack{
|\beta|=2} }
\|\langle ct-r \rangle
\partial^\beta u(t,x)\|_{L^2_x(\mathbb{R}^3)}
\\
\lesssim 
\sum_{ \substack{ \mu+|\alpha|\le 1} }
\|
L^\mu Z^\alpha u'(t,\cdot)\|_{L^2(\mathbb{R}^3)}
+
\|(t+r) (\partial_t^2-c^2\Delta)
u(t,x)\|_{L^2_x(\mathbb{R}^3)}
\end{multline}
by Sideris and Tu \cite[Lemma 7.1]{Sideris-Tu-2001-SIAM}.
See also Sogge \cite[p74, Lemma 5.3]{Sogge-2008-InterPress}.
\ebox

%
\subsection{Estimates for null conditions}

The null conditions are treated by the following estimates.

\begin{lemma}
\label{Lemma-NullConditions}
Let $c>0$. 
Let 
\beq
B(u,v)=\sum_{0\le j,k\le 3} B^{jk}\partial_ju\partial_kv,
\ \ \ \ 
Q(u,v)=\sum_{0\le j,k,l\le 3} Q^{jkl}\partial_ju\partial_k\partial_l v
\eeq 
satisfy the null conditions : 
\begin{equation}
\sum_{0\le j,k\le 3}
B^{jk}\xi_j \xi_k
=
\sum_{0\le j,k,l\le 3}
Q^{jkl}\xi_j \xi_k \xi_l=0
\ \ \ \ \mbox{for}\ \ \xi_0^2=c^2(\xi_1^2+\xi_2^2+\xi_3^2).
\end{equation}
Then the following inequalities hold for any $\alpha$ and functions $u$ and $v$.
\begin{equation}
(1) \ \ 
|\Gamma^\alpha B(u,v)|
\lesssim
\sum_{\substack{\beta+\gamma\le\alpha } }
\{
|\partialbar_c \Gamma^\beta u| |(\Gamma^\gamma v)'|
+
|(\Gamma^\beta u)'| |\partialbar_c \Gamma^\gamma v|
\},
\end{equation}
\begin{multline}
(2) \ \ 
|\Gamma^\alpha Q(u,v)|
\lesssim 
\sum_{\substack{\beta+\gamma\le\alpha } }
\Big\{
|\partialbar_c \Gamma^\beta u| |(\Gamma^\gamma v)''|
+
|(\Gamma^\beta u)'| |\partialbar_c(\Gamma^\gamma v)'|
\\
+
\frac{|(\Gamma^\beta u)'|
(
|(\Gamma^\gamma v)'|
+
|(\Gamma^\gamma v)''|
)
}{\angler} 
\Big\},
\end{multline}
where $\beta\le \alpha$ means any component of the multiindices satisfies the inequality.
\end{lemma}

\proof
(1) 
First, we consider the case $\alpha=0$.
Let $\omega_0=-c$ and $\omega=(\omega_1,\omega_2,\omega_3)\in \mathbb{S}^{2}$.
Since $\sum_{0\le j,k\le 3} B^{jk} \omega_j  \omega_k=0$ by the null condition, we have 
\begin{multline}
B(u,v)=
\sum_{0\le j,k\le 3} B^{jk} \partial_j u \partial_k v
-
\sum_{0\le j,k\le 3} B^{jk} \omega_j  \omega_k \partial_r u\partial_r v
\\
=
\sum_{0\le j,k\le 3} B^{jk} 
\{
(\partial_j  -\omega_j\partial_r)u\partial_k v
+
\omega_j\partial_ru(\partial_k-\omega_k\partial_r)v
\}.
\end{multline}
So that, we have 
$
|B(u,v)|\lesssim |\partialbar_c u||v'|+|u'||\partialbar_c v|
$.

For any $\alpha$, we have by the similar argument for Lemma 4.1 in \cite{Sideris-Tu-2001-SIAM} 
\begin{equation}
\Gamma^\alpha B(u,v)
= 
\sum_{\beta+\gamma\le \alpha}
B_{\beta,\gamma}(\Gamma^\beta u, \Gamma^\gamma v),
\end{equation}
where $\{B_{\beta,\gamma} \}_{\beta+\gamma\le \alpha}$ are quadratic nonlinear terms which satisfy the null conditions.
The required estimate follows from the above two results.

(2) By the same argument for (1), we have 
\begin{equation}
Q(u,v)=\sum_{0\le j,k,l\le 3} Q^{jkl}
\{
\partialbar_{cj}u\partial_k\partial_lv
+
\omega_j\partial_ru
(\partialbar_{ck}\partial_l v+\omega_k\partial_r\partialbar_{cl} v)
\}.
\end{equation}
So that, we have 
\begin{equation}
|Q(u,v)|\lesssim |\partialbar_c u||v''|+ |u'||\partialbar_c v'|+\frac{|u'||v'|}{r},
\end{equation}
where we have used 
$|\partial_r \partialbar_{cl} v|\lesssim |\partialbar_{cl} \nabla v|+|v'|/r$ 
for the last inequality.
Since $|Q(u,v)|\lesssim |u'||v''|$ for $r<1$, we have the required inequality for the case $\alpha=0$.
The case $\alpha\neq0$ also follows from 
\begin{equation}
\Gamma^\alpha Q(u,v)
= 
\sum_{\beta+\gamma\le\alpha}
Q_{\beta,\gamma}(\Gamma^\beta u, \Gamma^\gamma v),
\end{equation}
where $\{Q_{\beta,\gamma}\}_{\beta+\gamma\le \alpha} $ are quadratic nonlinear terms which satisfy the null conditions.
\newline \ \ \ebox

%

\newsection{Continuity argument to prove Theorem \ref{Theorem} }
We prepare the following proposition to prove Theorem \ref{Theorem}.

\begin{proposition}
\label{Prop}
Let $M_0$ and $N$ be positive numbers which satisfy 
 $M_0+5\le N\le 2M_0-2$.
For example, we are able to take $M_0=7$ and $N=12$.
We put 
\beq
\varepsilon:=
\sum_{|\alpha|\le N}
\|\anglex^{|\alpha|} \partial_x^\alpha \nabla f(x)\|_{L^2(\whole)} 
+
\sum_{|\alpha|\le N}
\|\anglex^{|\alpha|}\partial_x^{\alpha} g(x)\|_{L^2(\whole)}.
\eeq
Let $T>0$ and $A_0>0$. 
Let $u\in C^\infty([0,T)\times \mathbb{R}^3)$ be the local solution of \eqref{P-Whole}.
We assume 
\begin{equation}
\label{Ineq-Assumption}
\sum_{\substack{|\alpha|\le M_0 \\ 1\le I\le D} }
\sup_{\substack{0\le t< T \\ x\in \mathbb{R}^3} }
\angler^{1/2}\angletpr^{1/2}\anglecItmr^{1/2}
|\Gamma^\alpha u_I'(t,x)|\le A_0 \varepsilon.
\end{equation}
Then 
there exist constants $C_0>0$, which is independent of $A_0$, and $C>0$, which is dependent on $A_0$, such that the following estimates hold.
\beq
\label{Energy-High-R3}
(1)\ \ \ \ 
\sum_{|\alpha|\le M} \|\Gamma^\alpha u'(t,\cdot)\|_{L^2(\whole)}
\le C\varepsilon (1+t)^{C\varepsilon}
\ \ \ \ \mbox{for}\ \ 0\le t< T,\ 0\le M\le N.
\eeq
\begin{multline}
\label{Ineq-Linear-M2}
(2)\ \ \ \ 
\sum_{|\alpha|\le M-1}
(\log(e+t))^{-1/2}\|\anglex^{-1/2}\Gamma^\alpha u'(\cdot,x)\|_{L^2(S_t)}
\\
+
\sum_{\substack{|\alpha|\le M-1 \\ 1\le I\le D} }
(\log(e+t))^{-1/2}
\|\langle c_Is-r\rangle^{-1/2} \partialbar_{c_I} \Gamma^\alpha u_I(s,x)\|
_{L^2(S_t)} 
\\
\le
C\varepsilon (1+t)^{C\varepsilon}
\ \ \mbox{for}\ \ 0\le t< T,\ 0\le M\le N.
\end{multline}
\begin{multline}
\label{Ineq-Assumption-High}
(3)\ \ \ \ 
\sum_{\substack{|\alpha|\le M_0+3 \\ 1\le I\le D} } 
\sup_{x\in \whole}
\langle r\rangle^{1/2} \langle r+t \rangle^{1/2}
\anglecItmr^{1/2}
|\Gamma^\alpha u_I'(t, x)|
\\
\le
C\varepsilon (1+t)^{C\varepsilon}
\ \ \ \ \mbox{for}\ \ 0\le t< T.
\end{multline}
\beq
\label{Ineq-LowEnergy}
\ \hspace{-3.5cm}
(4)\ \ \ \ 
\sum_{|\alpha|\le M_0+2 }
\|\Gamma^\alpha u'\|_{L^\infty((0,T),L^2(\whole))}
\le C_0\varepsilon + C\varepsilon^{3/2}.
\eeq
\begin{multline}
(5)\ \ \ \ 
\sum_{\substack{|\alpha|\le M_0 \\ 1\le I\le D} } 
\sup_{\substack{0\le t< T \\ x\in \whole } }
\langle r\rangle^{1/2} \langle t+r \rangle^{1/2}
\anglecItmr^{1/2}
|\Gamma^\alpha u_I'(t, x)|
\le
C_0\varepsilon + C\varepsilon^{3/2}.
\end{multline}
\end{proposition}


\subsection{Proof of Theorem \ref{Theorem}}
Here we prove Theorem \ref{Theorem}.
We use the continuity argument which shows that the local in time solution $u$ does not blow up if its initial data is sufficiently small.
Since the constant $C_0$ is independent of $A_0$ in 
Proposition \ref{Prop}, 
we put $A_0=4C_0$ and take $\varepsilon$ sufficiently small such that 
$C\varepsilon^{3/2}\le C_0\varepsilon$.
Then the right hand side of (5) is bounded by 
$A_0\varepsilon/2$, which 
shows the local in time solution $u$ does not blow up, 
namely the solution exists globally in time.
\ebox


\subsection{Proof of Proposition \ref{Prop}}
First, we remark that under the assumption 
\eqref{Ineq-Assumption}, we have  
\begin{equation}
\label{Ineq-Assumption-Es}
\sum_{|\alpha|\le M/2+1}
\sup_{\substack{0\le t< T \\ x\in \whole} }
\angletpr
|\Gamma^\alpha u'(t,x)|\le C \varepsilon
\end{equation}
for some constant $C>0$ since $M/2+1\le M_0$.

(1)
For $1\le I,K\le D$ and $0\le k,l\le 3$, we put 
\beq
\gamma_I^{Kkl}
:=
\sum_{ \substack{1\le J\le D \\ 0\le j\le 3} }
Q_I^{JKjkl}\partial_ju_J.
\eeq
For any $\alpha$ with $|\alpha|\le M$, we use \eqref{Eqn-Divergence} and its integration to have 
\begin{multline}
\partial_t\int_{ \mathbb{R}^3 } e_0(\Gamma^\alpha u) dx
\lesssim
\sum_{1\le I\le D}
\|\square_{ \gamma_I } \Gamma^\alpha u_I\|_{L^2(\whole)}
\|(\Gamma^\alpha u_I)'\|_{ L^2(\mathbb{R}^3) }
\\
+
\sum_{ \substack{1\le I,K\le D \\ 0\le k,l\le 3} }
\|\partial_{t,x}\gamma_I^{Kkl}\|_{ L^\infty(\whole) }
\|(\Gamma^\alpha u_I)'\|_{ L^2(\mathbb{R}^3) }
\|(\Gamma^\alpha u_K)'\|_{ L^2(\mathbb{R}^3) }
. 
\end{multline}
To bound the right hand side, we use 
\begin{multline}
\label{Ineq-Commutator}
\sum_{\substack{1\le I\le D \\  |\alpha|\le M} }
\|\square_{ \gamma_I } \Gamma^\alpha u_I\|_{L^2(\whole)}
\lesssim 
\sum_{\substack{1\le I\le D \\  |\alpha|\le M} }
\| \Gamma^\alpha \square_{ \gamma_I }u_I\|_{L^2(\whole)}
+
\sum_{\substack{1\le I\le D \\  |\alpha|\le M-1} }
\| \Gamma^\alpha \square_{ c_I }u_I\|_{L^2(\whole)}
\\
+
\sum_{\substack{|\alpha|+|\beta|\le M \\ |\beta|\le M-1} }
\sum_{\substack{1\le I,K\le D \\ 0\le k,l\le 3} }
\|(\Gamma^\alpha \gamma_I^{Kkl}) \Gamma^\beta u_I''\|_{L^2(\whole)}
\lesssim
\frac{\varepsilon}{1+t} \sum_{|\alpha|\le M}\|\Gamma^\alpha u'\|_{L^2(\whole)},
\end{multline}
where we have used \eqref{Ineq-Assumption-Es} for the last inequality.
Since 
$\sum_{|\alpha|\le M}\|\Gamma^\alpha u'\|_{L^2(\whole)}$ 
is equivalent to  
$\sum_{|\alpha|\le M}\{\int_\whole e_0(\Gamma^\alpha u)dx\}^{1/2}$ 
for small $\varepsilon$, 
we obtain 
\beq
\partial_t 
\left\{
\sum_{\substack{|\alpha|\le M} }
\int_{ \mathbb{R}^3 } e_0(\Gamma^\alpha u) dx\right\}^{1/2}
\lesssim 
\frac{\varepsilon}{1+t} 
\left\{
\sum_{|\alpha|\le M} 
\int_{ \mathbb{R}^3 } e_0(\Gamma^\alpha u) dx\right\}^{1/2},
\eeq
which leads to the required inequality by the Gronwall inequality. 

%

(2)
By Lemma \ref{Lemma-EnergyKSSTangential}, the left hand side of \eqref{Ineq-Linear-M2} is bounded by  
\begin{equation}
C_0\varepsilon
+
C\sum_{\substack{1\le I\le D \\ |\alpha|\le M-1} }
\|\Gamma^\alpha \square_{c_I} u_I\|_{L^1((0,t),L^2(\mathbb{R}^3))}.
\end{equation}
For the last term, we use 
\begin{equation}
\label{Ineq-Nonlinear-Base}
\sum_{\substack{1\le I\le D \\ |\alpha|\le M-1 } }
|\Gamma^\alpha \square_{c_I} u_I|
\lesssim
\sum_{|\beta|\le M/2}|\Gamma^\beta u'|
\sum_{|\alpha|\le M}|\Gamma^\alpha u'|
\lesssim 
\frac{\varepsilon}{1+s} 
\sum_{|\alpha|\le M}|\Gamma^\alpha u'|,
\end{equation}
where we have used 
\eqref{Ineq-Assumption-Es} 
for the last inequality, and (1) to obtain 
\begin{multline}
\sum_{\substack{1\le I\le D \\ |\alpha|\le M-1} }
\|\Gamma^\alpha \square_{\gamma_I} u_I\|_{L^1((0,t),L^2(\mathbb{R}^3))}
\lesssim
\int_0^t \frac{\varepsilon}{1+s} \sum_{|\alpha|\le M}\|\Gamma^\alpha u'(s, \cdot)\|_{L^2(\mathbb{R}^3)} dt
\\
\lesssim 
\varepsilon (1+t)^{C\varepsilon}.
\end{multline}
Therefore, we obtain the required inequality.

%

(3) 
Let $M=M_0+3$.
By Lemma \ref{Lemma-KS}, the left hand side of \eqref{Ineq-Assumption-High} is bounded by  
\begin{multline}
C_0\sum_{|\alpha|\le M+2} 
\|\Gamma^\alpha u'(t,\cdot)\|_{L^2(\mathbb{R}^3)}
+
C_0\sum_{ \substack{|\alpha|\le M+1 } } 
\|\langle t+|y|\rangle \Gamma^\alpha \square_{c_I} u_I(t,y)\|_{L^2_y(\mathbb{R}^3)}.
\end{multline}
The last term is bounded by 
$
C\varepsilon 
\sum_{|\alpha|\le M+2} 
\|\Gamma^\alpha u'(t,\cdot)\|_{L^2(\mathbb{R}^3)}
$
by \eqref{Ineq-Assumption-Es} since $(M+2)/2\le M_0$.
This shows the required inequality by (1).

%

(4) 
Let $M=M_0+2$.
By the standard energy estimates, we have 
\begin{multline}
\sum_{|\alpha|\le M}
\|\Gamma^\alpha u_I'(t,\cdot)\|_{L^2(\whole)}^2
\le 
C_0\sum_{|\alpha|\le M}
\|\Gamma^\alpha u_I'(0,\cdot)\|_{L^2(\whole)}^2
\\
+
C_0\sum_{|\alpha|\le M}
\int_0^t\int_{\whole} 
|\partial_t \Gamma^\alpha u_I \square_{c_I}\Gamma^\alpha u_I|dxds
=: A_1+A_2
\end{multline}
for $1\le I\le D$.
We have $A_1\le (C_0\varepsilon)^2$ for some $C_0>0$ which is independent of $A_0$.
By Lemma \ref{Lemma-NullConditions}, $A_2$ is bounded as  
\begin{multline}
A_2\lesssim
\sum_{\substack{1\le J,K\le D \\ (c_J,c_K)=(c_I,c_I) } }
\int_0^t\int_\whole
\sum_{|\alpha|\le M} |\Gamma^\alpha u_I'|
\sum_{|\alpha|\le M+1} |\partialbar_{c_J}\Gamma^\alpha u_J|
\sum_{|\alpha|\le M+1} |\Gamma^\alpha u_K'|
dxds
\\
+
\sum_{\substack{1\le J,K\le D \\ (c_J,c_K)=(c_I,c_I) } }
\int_0^t\int_\whole
\sum_{|\alpha|\le M} |\Gamma^\alpha u_I'|
\sum_{|\alpha|\le M+1} |\Gamma^\alpha u_J'|
\sum_{|\alpha|\le M+1} |\Gamma^\alpha u_K'|
\frac{dx}{\angler}ds
\\
+
\sum_{\substack{1\le J,K\le D \\ (c_J,c_K)\neq(c_I,c_I) } }
\int_0^t\int_\whole
\sum_{|\alpha|\le M} |\Gamma^\alpha u_I'|
\sum_{|\alpha|\le M+1} |\Gamma^\alpha u_J'|
\sum_{|\alpha|\le M+1} |\Gamma^\alpha u_K'|
dxds
\\
=: A_3+A_4+A_5.
\end{multline}
We use \eqref{Ineq-Assumption-High} to have 
\begin{multline}
A_3
\lesssim 
\varepsilon
\sum_{\substack{1\le J,K\le D \\ (c_J,c_K)=(c_I,c_I) } }
\sum_{|\alpha|\le M+1}
\|\langle c_Is-r \rangle^{-1/2}\langle s \rangle^{-\delta} \partialbar_{c_J} \Gamma^\alpha u_J\|_{L^2((0,t)\times\whole)}
\\
\cdot
\sum_{|\alpha|\le M+1} \|\angler^{-1/2}\langle s \rangle^{-\delta} \Gamma^\alpha u_K'\|_{L^2((0,t)\times\whole)}
\lesssim 
\varepsilon^3,
\end{multline}
where $\delta>0$ is a sufficiently small number and we have used 
\eqref{Ineq-Linear-M2} to obtain the last inequality.
Similarly, we have 
\begin{multline}
A_4
\lesssim 
\varepsilon
\sum_{\substack{1\le J,K\le D \\ (c_J,c_K)=(c_I,c_I) } }
\sum_{|\alpha|\le M+1}
\|\angler^{-1/2}\langle s \rangle^{-\delta} \Gamma^\alpha u_J'\|_{L^2(S_t)}
\\
\cdot
\sum_{|\alpha|\le M+1} \|\angler^{-1/2}\langle s \rangle^{-\delta} \Gamma^\alpha u_K'\|_{L^2(S_t)}
\lesssim 
\varepsilon^3.
\end{multline}
To bound $A_5$, we consider the conic neighborhood defined by 
\beq
\Lambda_I:=\{(s,x)\ :\ 0\le s\le t,\ |c_Is-r|\le c_0s/10\},
\ \ \ \ 
c_0:=\max_{1\le J\le D} c_J
\eeq
for $1\le I\le D$.
We note that 
\beq
\sum_{|\alpha|\le M+1} |\Gamma^\alpha u_I'(s,x)|
\lesssim 
\varepsilon \angler^{-1/2}{\langle s+r \rangle}^{-1+C\varepsilon}
\eeq
on $\whole\backslash \Lambda_I$ by \eqref{Ineq-Assumption-High} for $1\le I\le D$.
So that, we have  
\begin{multline}
\int_0^t\int_{\whole\backslash\Lambda_I}
\sum_{|\alpha|\le M+1} |\Gamma^\alpha u_I'|
\sum_{|\alpha|\le M+1} |\Gamma^\alpha u_J'|
\sum_{|\alpha|\le M+1} |\Gamma^\alpha u_K'|
dxds
\\
\lesssim 
\varepsilon 
\sum_{|\alpha|\le M+1}
\|\angler^{-1/2}\angles^{-\delta} \Gamma^\alpha u_J'\|_{L^2(S_t)}
\\
\cdot
\sum_{|\alpha|\le M+1}
\|\angler^{-1/2}\angles^{-\delta} \Gamma^\alpha u_K'\|_{L^2(S_t)}
\lesssim 
\varepsilon^3,
\end{multline}
where $\delta>0$ is sufficiently small.
Since 
$\whole = 
(\whole\backslash \Lambda_I)
\cup 
(\whole\backslash \Lambda_J)
\cup 
(\whole\backslash \Lambda_K)$ 
by $(c_J,c_K)\neq(c_I,c_I)$, 
we obtain $A_5\lesssim\varepsilon^3$ and the required inequality.

%
(5) 
The left hand side of (5) is bounded by $(C_0+C\varepsilon)\sum_{|\alpha|\le M_0+2}\|\Gamma^\alpha u'(t,\cdot)\|_{L^2(\mathbb{R}^3)}$ by the similar argument for the proof of (3). The required inequality follows from (4).
\ebox


\newsection{Several estimates to prove Theorem \ref{Theorem-Ext} }
In this section, we prepare several estimates to prove Theorem \ref{Theorem-Ext} in the next section.
Let $c>0$, $0<T\le \infty$.
We show the estimates for the solution of the scalar-valued problem 
\begin{equation}
\label{Eqn-Cauchy-F}
\left\{
\begin{array}{l}
\square_c u=F \ \ \mbox{for}\ \ (t,x)\in [0,T)\times \ext \\
u(t,\cdot)|_{\partial \mathcal{K}} =0 \ \ \mbox{for}\ \ t\in [0,T) \\
u(0,\cdot)=f(\cdot),\ \ \partial_tu(0,\cdot)=g(\cdot),
\end{array}
\right.
\end{equation}
where $\square_c:=\partial_t^2-c^2\Delta$.
Let $\zeta \in C_0^\infty(\whole)$ be a function which satisfies 
$0\le \zeta \le1$, 
$\zeta (x)=1$ for $|x|\le 3$, and $\zeta (x)=0$ for $|x|\ge4$.
Regarding $(1-\zeta) F$, $(1-\zeta) f$, $(1-\zeta) g$ as  functions 
on $\whole$ by zero-extension, let $v$ be the solution of 
\begin{equation}
\label{Eqn-Cauchy-v}
\left\{
\begin{array}{l}
\square_c v=(1-\zeta)F \ \ \mbox{for}\ \ (t,x)\in [0,T)\times \whole \\
v(0,\cdot)=((1-\zeta) f)(\cdot),\ \ \partial_tv(0,\cdot)=((1-\zeta) g)(\cdot).
\end{array}
\right.
\end{equation}

\subsection{Estimates for boundary terms}

We use the following estimates to bound the terms from the commutator estimates.

\begin{lemma}
\label{Lemma-BoundaryExponential}
The solution $u$ of \eqref{Eqn-Cauchy-F} satisfies the following estimate.
\begin{multline}
\label{Ineq-BoundaryExponential}
\sum_{\substack{\mu+|\alpha|\le M \\ \mu\le \mu_0} }
\|L^\mu \partial^\alpha u'(t,\cdot)\|_{L^2(\{x\in \ext: |x|<2\})}
\lesssim
e^{-at/2}
\sum_{\substack{\mu+|\alpha|\le M+2 \\ \mu\le \mu_0} }
\|L^\mu \partial^\alpha u(0,\cdot)\|_{L^2(\ext)}
\\
+
\int_0^t e^{-a(t-s)/2}
\sum_{\substack{\mu+|\alpha|\le M+1 \\ \mu\le \mu_0} }
\|L^\mu \partial^\alpha \square_c u(s,\cdot)\|_{L^2(\{x\in \ext: |x|<4\})}
ds
\\
+
\sum_{\substack{\mu+|\alpha|\le M-1 \\ \mu\le \mu_0} }
\|L^\mu \partial^\alpha \square_c u(t,\cdot)\|_{L^2(\{x\in \ext: |x|<4\})}
\\
+
\int_0^t e^{-a(t-s)/2}
\sum_{\substack{\mu+|\alpha|\le M+1 \\ \mu\le \mu_0, \ |\beta|\le1} }
\|L^\mu \partial^\alpha \partial^\beta v(s,\cdot)\|_{L^2(\{x\in \whole: |x|<3\})}
ds
\\
+
\sum_{\substack{\mu+|\alpha|\le M-1 \\ \mu\le \mu_0,\ |\beta|\le 1} }
\|L^\mu \partial^\alpha \partial^\beta v(t,\cdot)\|_{L^2(\{x\in \whole : |x|<3\})}
\ \ \ \ \mbox{for}\ \ 0\le t<T.
\end{multline}
\end{lemma}
\proof
For simplicity, we only consider the case $c=1$. 
Let $u_1$ and $u_2$ be the solutions of the boundary value problems 
\beq
\left\{
\begin{array}{l}
\square u_1=\zeta \square u,\ \ 
u_1|_{\partial \mathcal{K}}=0,\\ 
u_1(0,\cdot)=\zeta u(0,\cdot), \ \ 
\partial_t u_1(0,\cdot)=\zeta \partial_tu(0,\cdot),
\end{array}
\right.
\eeq
\beq
\left\{
\begin{array}{l}
\square u_2=(1-\zeta) \square u,\ \ 
u_2|_{\partial \mathcal{K}}=0,\\ 
u_2(0,\cdot)=(1-\zeta) u(0,\cdot), \ \ 
\partial_t u_2(0,\cdot)=(1-\zeta) \partial_tu(0,\cdot).
\end{array}
\right.
\eeq
Then we have $u=u_1+u_2$.
By the local energy decay estimates \eqref{LED}, we have 
\begin{multline}
\sum_{\substack{\mu+|\alpha|\le M \\ \mu\le \mu_0} }
\|L^\mu \partial^\alpha u_1'(t,\cdot)\|_{L^2(\{x\in \ext:|x|<2\})}
\lesssim 
e^{-at/2} 
\sum_{\substack{\mu+|\alpha|\le M+2 \\ \mu\le \mu_0} }
\|L^\mu \partial^\alpha u(0,\cdot)\|_{L^2(\ext)}
\\
+
\int_0^t e^{-a(t-s)/2}
\sum_{\substack{\mu+|\alpha|\le M+1 \\ \mu\le \mu_0} }
\|L^\mu \partial^\alpha \square u(s,\cdot)\|_{L^2( \{x\in \ext:|x|<4\} )}
ds
\\
+
\sum_{\substack{\mu+|\alpha|\le M-1 \\ \mu\le \mu_0} }
\|L^\mu \partial^\alpha \square u(t,\cdot)\|_{L^2( \{x\in \ext:|x|<4\} )}.
\end{multline}
Let $u_3:=v|_{\ext}$ and  $u_4:=u_2-u_3$.
Then $\square u_4=u_4(0,\cdot)=\partial_t u_4(0,\cdot)=0$.
Let $\rho\in C_0^\infty(\whole)$ be a function with $\rho(x)=1$ for $|x|\le 2$, and $\rho(x)=0$ for $|x|\ge3$.
We put $\tilde{u_2}:=\rho u_3+u_4$.
Then $\tilde{u_2}=u_2$ for $|x|\le 2$, and 
\begin{multline}
\square \tilde{u_2}=-2\nabla\rho\cdot \nabla u_3-(\Delta \rho)u_3,\ \ 
\tilde{u_2}(0,\cdot)=\rho u_3(0,\cdot),\ \ 
\partial_t \tilde{u_2}(0,\cdot)=\rho \partial_t u_3(0,\cdot)
\end{multline}
So that, by \eqref{LED}, we have 
\begin{multline}
\sum_{\substack{\mu+|\alpha|\le M \\ \mu\le \mu_0} }
\|L^\mu \partial^\alpha u_2'(t,\cdot)\|_{L^2( \{x\in \ext : |x|<2\} )}
=
\sum_{\substack{\mu+|\alpha|\le M \\ \mu\le \mu_0} }
\|L^\mu \partial^\alpha \tilde{u_2}'(t,\cdot)\|_{L^2( \{x\in \ext : |x|<2\} )}
\\
\lesssim 
e^{-at/2}
\sum_{\substack{\mu+|\alpha|\le M+2 \\ \mu\le \mu_0} }
\|L^\mu \partial^\alpha u(0,\cdot)\|_{L^2(\ext)}
\\
+
\int_0^t e^{-a(t-s)/2}
\sum_{\substack{\mu+|\alpha|\le M+1 \\ \mu\le \mu_0,\ |\beta|\le 1} }
\|L^\mu \partial^\alpha \partial^\beta v(s,\cdot)\|_{L^2(\{x\in \whole : |x|<3\})}
ds
\\
+
\sum_{\substack{\mu+|\alpha|\le M-1 \\ \mu\le \mu_0,\ |\beta|\le 1} }
\|L^\mu \partial^\alpha \partial^\beta v(t,\cdot)\|_{L^2(\{x\in \whole : |x|<3\})}.
\end{multline}
Therefore, we obtain the required estimate.
\ebox

\begin{lemma}
\label{Lemma-Boundary-t}
The solution $u$ of \eqref{Eqn-Cauchy-F} satisfies the following estimate.
\begin{multline}
(1)\ \ \sum_{\substack{\mu+|\alpha|\le M \\ \mu\le \mu_0} }
(1+t)\|L^\mu \partial^\alpha u'(t,\cdot)\|_{L^2(\{x\in \ext : |x|<2\})}
\lesssim
\sum_{\substack{\mu+|\alpha|\le M+2 \\ \mu\le \mu_0} }
\|L^\mu \partial^\alpha u(0,\cdot)\|_{L^2(\ext)}
\\
+
\sum_{\substack{\mu+|\alpha|\le M+1 \\ \mu\le \mu_0} }
\sup_{0\le s\le t}
(1+s)\|L^\mu Z^\alpha \square_c u(s,\cdot)\|_{L^2(\{x\in \ext: |x|<4\})}
\\
+
\sum_{\substack{\mu+|\alpha|\le M+1 \\ \mu\le \mu_0,\ |\beta|\le 1} }
\sup_{0\le s\le t}
(1+s)\|L^\mu \partial^\alpha \partial^\beta v(s,\cdot)\|_{L^2( \{x\in \whole : |x|<3\} )}
\ \ \mbox{for}\ \ 0\le t<T.
\end{multline}
\begin{multline}
(2)\ \ \sum_{\substack{\mu+|\alpha|\le M \\ \mu\le \mu_0} }
\|L^\mu \partial^\alpha u'\|_{L^2((0,T)\times\{x\in \ext : |x|<2\})}
\lesssim
\sum_{\substack{\mu+|\alpha|\le M+2 \\ \mu\le \mu_0} }
\|L^\mu \partial^\alpha u(0,\cdot)\|_{L^2(\ext)}
\\
+
\sum_{\substack{\mu+|\alpha|\le M+1 \\ \mu\le \mu_0} }
\|L^\mu \partial^\alpha \square_c u\|_{L^1((0,T), L^2(\{x\in \ext: |x|<4\}))}
\\
+
\sum_{\substack{\mu+|\alpha|\le M-1 \\ \mu\le \mu_0} }
\|L^\mu \partial^\alpha \square_c u\|_{L^2( (0,T)\times\{x\in \ext: |x|<4\} )}
\\
+
\sum_{\substack{\mu+|\alpha|\le M+1 \\ \mu\le \min\{\mu_0,M\},\ |\beta|\le1 } }
\|L^\mu \partial^\alpha \partial^\beta v\|_{L^2( (0,T)\times \{x\in \whole : |x|<3\} )}
\end{multline}

\end{lemma}

\proof
The results follow from Lemma \ref{Lemma-BoundaryExponential}, the boundedness of  $te^{-at/2}$ for (1), and the Young inequality for the time variable for (2).
\ebox

\vspace{10pt}

The following result has been partially shown in 
\cite[Lemma 2.9]{Metcalfe-Sogge-2005-InventMath}
for vanishing Cauchy data.
We generalize it in complete form.

\begin{lemma}
\label{Lemma-BoundaryTerm}
The solution $u$ of \eqref{Eqn-Cauchy-F} satisfies 
\begin{multline}
\sum_{\substack{\mu+|\alpha|\le M \\ \mu\le \mu_0}}
\int_0^t \|L^\mu \partial^\alpha u'(s,x)\|_{L^2( \{ x\in \ext : |x|<2 \} )}ds
\lesssim
\sum_{\substack{\mu+|\alpha|\le M+2 \\ \mu\le \mu_0} }
\|\anglex L^\mu\partial^\alpha u(0,x)\|_{L^2_x(\ext)}
\\
+
\sum_{\substack{\mu+|\alpha|\le M+1 \\ \mu\le \mu_0}}
\int_0^t \|L^\mu\partial^\alpha \square_cu(\tau,\cdot)\|_{L^2(\{x\in\ext:|x|<4\})} d\tau 
\\
+
\sum_{\substack{\mu+|\alpha|\le M+1 \\ \mu\le \mu_0}}
\int_0^t\int_0^s \|L^\mu\partial^\alpha \square_c u(\tau,\cdot)\|_{L^2(\{x\in\ext : ||x|-c(s-\tau)|<4\} )} d\tau ds.
\end{multline}
\end{lemma}

\proof
The required result easily follows from the integration by $t$ of 
\eqref{Ineq-BoundaryExponential} and the Young inequality if we show 
\begin{multline}
\label{Proof-BT-Ineq}
\sum_{\substack{\mu+|\alpha|\le M+1 \\ \mu\le \mu_0,\ |\beta|\le 1} }
\int_0^t \|L^\mu \partial^\alpha \partial^\beta v(s,\cdot)\|_{L^2({x\in \mathbb{R}^3 : |x|<3})} ds
\\
\lesssim
\sum_{\substack{\mu+|\alpha|\le M+2 \\ \mu\le \mu_0} }
\|\anglex L^\mu \partial^\alpha u(0,x)\|_{L^2_x(\ext)}
\\
+
\sum_{\substack{\mu+|\alpha|\le M+1 \\ \mu\le \mu_0} }
\int_0^t\int_0^s
\|L^\mu \partial^\alpha \square_cu(\tau,\cdot)\|_{L^2(\{x\in \ext : ||x|-c(s-\tau)|<4\})} d\tau ds.
\end{multline}
To show this inequality, we prepare the following claim.

\begin{claim}
For any given functions $w_0$, $w_1$ and $G$, the solution $w$ of the Cauchy problem 
\beq
\left\{
\begin{array}{l}
\square_c w =G \ \ \mbox{on}\ [0,\infty)\times \mathbb{R}^3 \\
w(0,\cdot)=w_0(\cdot),\ \partial_t w(0,\cdot)=w_1(\cdot)
\end{array}
\right.
\eeq
satisfies 
\begin{multline}
\sum_{|\alpha|\le 1} \|\partial^\alpha w(t,\cdot)\|_{L^2(\{x\in \mathbb{R}^3 : |x|<3\})}
\lesssim
\sum_{|\alpha|\le 1} \|\partial^\alpha w(0,\cdot)\|_{L^2(\{x\in \mathbb{R}^3 : ||x|-ct|<4\})}
\\
+
\int_0^t \|\square_c w(s,\cdot)\|_{L^2(\{x\in \mathbb{R}^3 : ||x|-c(t-s)|<4\})} ds 
\end{multline}
for $t\ge0$.
\end{claim}

\proof 
We note that $w$ is written as 
\beq
w(t,\cdot)=\partial_t K(t) w_0(\cdot)+K(t)w_1(\cdot)+\int_0^t K(t-s)G(s,\cdot)ds,
\eeq
where $K(t):=\sin ct\sqrt{-\Delta}/c\sqrt{-\Delta}$, $\partial_t K(t):=\cos ct\sqrt{-\Delta}$.
Let $\chi \in C_0^\infty(\mathbb{R})$ be a function with $\chi(s)=1$ for $-3\le s\le 3$ 
and $\chi(s)=0$ for $|s|\ge4$.
For any fixed $t\ge0$, by the Huygens principle, we have 
\begin{multline}
w(t,x)=
\partial_tK(t) \left( \chi(|\cdot|-ct)w_0(\cdot) \right) (x)
+
K(t) \left( \chi(|\cdot|-ct)w_1(\cdot) \right) (x)
\\
+
\int_0^t K(t-s) \left( \chi(|\cdot|-c(t-s))G(s,\cdot) \right) (x) ds
\end{multline}
for $x$ with $|x|<3$.
By the isometry of $e^{it\sqrt{-\Delta}}$ and $\nabla/\sqrt{-\Delta}$ on $L^2(\mathbb{R}^3)$ and the embedding 
$\dot{H}^1(\mathbb{R}^3) \hookrightarrow L^6(\mathbb{R}^3) \hookrightarrow L^2(\{x\in \mathbb{R}^3: |x|<3\})$, we have 
\begin{multline}
\label{Claim-BT-Ineq1}
\|w(t,\cdot)\|_{ L^2(\{x\in \mathbb{R}^3 : |x|<3\}) }
\lesssim 
\|w_0\|_{ L^2(\{x\in \mathbb{R}^3 : ||x|-ct|<4\}) }
\\
+
\|w_1\|_{ L^2(\{x\in \mathbb{R}^3 : ||x|-ct|<4\}) }
+
\int_0^t \|G(s,\cdot)\|_{ L^2(\{x\in \mathbb{R}^3 : ||x|-c(t-s)|<4\}) } ds
\end{multline}
and 
\begin{multline}
\label{Claim-BT-Ineq2}
\|\nabla w(t,\cdot)\|_{ L^2(\{x\in \mathbb{R}^3 : |x|<3\}) }
\lesssim 
\sum_{|\alpha|\le 1}
\|\partial^\alpha_x w_0\|_{ L^2(\{x\in \mathbb{R}^3 : ||x|-ct|<4\}) }
\\
+
\|w_1\|_{ L^2(\{x\in \mathbb{R}^3 : ||x|-ct|<4\}) }
+
\int_0^t \|G(s,\cdot)\|_{ L^2(\{x\in \mathbb{R}^3 : ||x|-c(t-s)|<4\}) } ds.
\end{multline}
Since $\partial_tw$ is written as  
\beq
\partial_t w(t,\cdot)=c^2 \Delta K(t)w_0(\cdot)+\partial_t K(t)w_1(\cdot)+\int_0^t \partial_t K(t-s)G(s,\cdot) ds,
\eeq
we have by the Huygens principle   
\begin{multline}
\partial_t w(t,x)=
c^2 \Delta K(t) \left( \chi(|\cdot|-ct)w_0(\cdot) \right) (x)
+
\partial_t K(t) \left( \chi(|\cdot|-ct)w_1(\cdot) \right) (x)
\\
+
\int_0^t \partial_tK(t-s) \left( \chi(|\cdot|-c(t-s))G(s,\cdot) \right) (x) ds
\end{multline}
for $x$ with $|x|<3$.
So that, $\|\partial_t w(t,\cdot)\|_{L^2(\{x\in \mathbb{R}^3 : |x|<3\})}$ is bounded by the right hand side of \eqref{Claim-BT-Ineq2}.
And we obtain the required inequality.
\ebox

We apply the above claim to $v$ and integrate it by $t$ to obtain 
\begin{multline}
\sum_{|\beta|\le1} 
\int_0^t \|\partial^\beta v(s,\cdot)\|_{ L^2(\{x\in \mathbb{R}^3 : |x|<3\}) } ds 
\lesssim 
\sum_{|\beta|\le 1}
\int_0^t \|\partial^\beta v(0,\cdot)\|_{ L^2(\{x\in \mathbb{R}^3 : ||x|-cs|<4\}) } ds 
\\
+
\int_0^t \int_0^s \|\square_c v(\tau,\cdot)\|_{ L^2(\{x\in \mathbb{R}^3 : ||x|-c(s-\tau)|<4\}) } d\tau ds.
\end{multline}
Here, the first term in the right hand side is bounded by 
\begin{multline}
\int_0^t \frac{1}{\langle s \rangle} \sum_{|\beta|\le 1} 
\| \langle s \rangle\partial^\beta v(0,\cdot)\|_{L^2(\{x\in \mathbb{R}^3:||x|-cs|<4\})} ds 
\lesssim
\sum_{|\beta|\le 1} \|\langle x \rangle\partial^\beta v(0,x)\|_{L^2_x(\mathbb{R}^3)}.
\end{multline}
By the replacement of $v$ with $L^\mu \partial^\alpha v$ and the definition of $v$, 
we obtain \eqref{Proof-BT-Ineq}.
\ebox

%

\subsection{Weighted energy estimates}

The following weighted energy estimates are the exterior domain 
analog to 
Lemma \ref{Lemma-EnergyKSSTangential}.

\begin{lemma}
\label{Lemma-WeightedEnergy-Ext}
For any $M\ge0$ and $\mu_0\ge0$, the solution $u$ of \eqref{Eqn-Cauchy-F} satisfies 
\begin{multline}
\sum_{\substack{\mu+|\alpha|\le M \\ \mu\le \mu_0}}
\|L^\mu Z^\alpha u'\|_{L^\infty((0,T), L^2(\ext))}
\\
+\sum_{\substack{\mu+|\alpha|\le M \\ \mu\le \mu_0}}
(\log(e+T))^{-1/2}
\|\anglex^{-1/2}L^\mu Z^\alpha u'\|_{L^2((0,T)\times \ext)}
\\
+
\sum_{\substack{\mu+|\alpha|\le M \\ \mu\le \mu_0}}
(\log(e+T))^{-1/2}
\|\langle ct-r \rangle^{-1/2}
\partialbar_{c} L^\mu Z^\alpha u\|_{L^2((0,T)\times \ext)}
\\
\lesssim 
\sum_{\substack{\mu+|\alpha|\le M+2 \\ \mu\le \mu_0} }
\|L^\mu Z^\alpha u(0,\cdot)\|_{L^2(\ext)}
+
\sum_{\substack{\mu+|\alpha|\le M+1 \\ \mu\le \mu_0}}
\|L^\mu Z^\alpha \square_c u\|_{L^1((0,T), L^2(\ext))} 
\\
+
\sum_{\substack{\mu+|\alpha|\le M-1 \\ \mu\le \mu_0}}
\|L^\mu \partial^\alpha \square_c u\|_{ L^2((0,T)\times\{x\in\ext:|x|<4\} ) }.
\end{multline}
Here, the above estimate holds with all $Z$ replaced by $\partial$.
\end{lemma}

\proof
The inequality for the second term in the left hand side has been shown by Metcalfe and Sogge \cite[Proposition 2.6]{Metcalfe-Sogge-2005-InventMath}.
We show the inequality for the third term.
The inequality for the first term follows similarly.
We only show the case $c=1$ for simplicity.
First, we consider the boundaryless case.

\begin{claim} 
\label{Lemma-EnergyKSSTangential-Claim1}
Let us consider the Cauchy problem 
\begin{equation}
\left\{
\begin{array}{l}
\square u=F \ \ \ \ \mbox{on}\ \ [0,T)\times \mathbb{R}^3 \\
u(0,\cdot)=f(\cdot),\ \ \partial_t u(0,\cdot)=g(\cdot).
\end{array}
\right.
\end{equation}
Let $\supp F\subset \{(t,x): 0\le t< T,\ |x|\le 2\}$.
Then we have 
\begin{multline}
(\log(e+T))^{-1/2}
\|\angletmr^{-1/2}\partialbar u\|_{L^2((0,T)\times\mathbb{R}^3)}
\\
\lesssim 
\|\nabla f\|_{L^2(\mathbb{R}^3)}
+\|g\|_{L^2(\mathbb{R}^3)}
+\|\square u\|_{L^2((0,T)\times\mathbb{R}^3)}.
\end{multline}
\end{claim}

\proof
Let $\chi\ge0$ be a bump function such that $\chi(t)=1$ for $-1/4\le t\le 1/4$, $\chi(t)=0$ for $t\le -1$ or $t\ge1$, 
and $\sum_{j=0}^\infty \chi(t-j)=1$ for any $t\ge0$.
Let $\{u_j\}_{j\ge-1}$ be the solutions of 
\begin{equation}
\square u_{-1}=0,\ \ u_{-1}(0,\cdot)=f(\cdot),\ \ \partial_tu_{-1}(0,\cdot)=g(\cdot),
\end{equation}
\begin{equation}
\square u_{j}(t,\cdot)=\chi(t-j)F(t,\cdot),\ \ u_{j}(0,\cdot)=0,\ \ \partial_tu_{j}(0,\cdot)=0
\end{equation}
for $j\ge0$.
Then we have $u=u_{-1}+\sum_{j=0}^\infty u_j$.
By the Huygens principle, for any  $(t,x)\in [0,T)\times\mathbb{R}^3$, there exists $j(t,x)\ge-1$ such that 
\begin{equation}
u(t,x)=u_{-1}(t,x)+\sum_{|j-j(t,x)|\le 3} u_j(t,x).
\end{equation}
So that, we have 
$
|\partialbar u(t,x)|^2
\lesssim 
\sum_{j\ge-1} |\partialbar u_j(t,x)|^2,
$ 
and  
\begin{multline}
(\log(e+T))^{-1}
\|\angletmr^{-1/2}\partialbar u\|_{L^2((0,T)\times\mathbb{R}^3))}^2
\\
\lesssim
\sum_{j=-1}^\infty
(\log(e+T))^{-1}
\|\angletmr^{-1/2}\partialbar u_j\|_{L^2((0,T)\times\mathbb{R}^3))}^2
\\
\lesssim 
\|\nabla f\|_{L^2(\mathbb{R}^3)}^2+C\|g\|_{L^2(\mathbb{R}^3)}^2
+
\sum_{j\ge0}\|\chi(\cdot-j)F\|_{L^1((0,T),L^2(\mathbb{R}^3))}^2,
\end{multline}
where we have used \eqref{Ineq-EnergyKSSTangential} for the last inequality.
Since the last term is bounded by $\|F\|_{L^2((0,T)\times \mathbb{R}^3)}^2$ due to the H{\"o}lder inequality, we obtain the required result.
\ebox

%
\begin{claim} 
\label{Lemma-EnergyKSSTangential-Claim2}
Let us consider the problem 
\begin{equation}
\left\{
\begin{array}{l}
\square u=F \ \ \ \ \mbox{on}\ \ [0,T)\times \ext \\
u(0,\cdot)=f(\cdot),\ \ \partial_t u(0,\cdot)=g(\cdot) \\
u(t,\cdot)|_{\mathcal{K}} =0 \ \ \ \ \mbox{for }\ \ t\in [0,T).
\end{array}
\right.
\end{equation}
Then we have 
\begin{multline}
\sum_{\substack{\mu+|\alpha|\le M \\ \mu\le \mu_0} }
(\log(e+T))^{-1/2}
\|\angletmr^{-1/2}\partialbar L^\mu Z^\alpha u\|
_{L^2((0,T)\times\ext)}
\\
\lesssim 
\sum_{\substack{\mu+|\alpha|\le M+1 \\ \mu\le \mu_0} }
\|L^\mu Z^\alpha u(0,\cdot)\|_{L^2(\ext)}
+
\sum_{\substack{\mu+|\alpha|\le M \\ \mu\le \mu_0} }
\|L^\mu Z^\alpha \square u\|_{L^1((0,T), L^2(\ext))}
\\
+
\sum_{\substack{\mu+|\alpha|\le M \\ \mu\le \mu_0} }
\|L^\mu \partial^\alpha u'\|_{L^2((0,T)\times\{x\in\ext:|x|<2\})}.
\end{multline}
\end{claim}

\proof
When $|x|\le 2$, the left hand side is bounded by the last term.
We consider the case $|x|\ge2$.
Let $0\le\eta\le1$ be a smooth function such that $\eta(x)=0$ for $|x|\le 1$, $\eta(x)=1$ for $|x|\ge2$. 
Then $\eta u$ can be seen as a function on $\mathbb{R}^3$ and satisfies 
\begin{equation}
\square (\eta u)=\eta \square u-2\nabla \eta\cdot \nabla u-(\Delta \eta)\cdot u.
\end{equation}
Let $u_1$ and $u_2$ be the solutions of 
\begin{equation}
\square u_1=\eta \square u,\ \ u_1(0,\cdot)=(\eta u)(0,\cdot),\ \ \partial_tu_1(0,\cdot)=(\eta \partial_tu)(0,\cdot),
\end{equation}
\begin{equation}
\square u_2=-2\nabla \eta\cdot\nabla u-(\Delta\eta)u,\ \ u_2(0,\cdot)=0,\ \ \partial_tu_2(0,\cdot)=0
\end{equation}
on $[0,T)\times \mathbb{R}^3$.
Since $u=\eta u$ for $|x|\ge2$, and $\eta u=u_1+u_2$, we have 
\begin{multline}
\sum_{\substack{\mu+|\alpha|\le M \\ \mu\le \mu_0} }
(\log(e+T))^{-1/2}
\|\angletmr^{-1/2}\partialbar L^\mu Z^\alpha u\|
_{L^2((0,T)\times\{x\in \ext : |x|>2\})}
\\
\le 
\sum_{j=1,2}
\sum_{\substack{\mu+|\alpha|\le M \\ \mu\le \mu_0} }
(\log(e+T))^{-1/2}
\|\angletmr^{-1/2}\partialbar L^\mu Z^\alpha u_j\|
_{L^2((0,T)\times\whole)}
=:I_1+I_2.
\end{multline}
For $u_1$, we use Lemma \ref{Lemma-EnergyKSSTangential} to obtain 
\begin{multline}
I_1
\lesssim 
\sum_{\substack{\mu+|\alpha|\le M \\ \mu\le \mu_0} }
\|L^\mu Z^\alpha u'(0,\cdot)\|_{L^2(\ext)}
+
\sum_{\substack{\mu+|\alpha|\le M \\ \mu\le \mu_0} }
\|L^\mu Z^\alpha \square u\|_{L^1((0,T), L^2(\ext))}.
\end{multline}
For $u_2$, we use 
Claim \ref{Lemma-EnergyKSSTangential-Claim1} 
to obtain 
\begin{multline}
I_2
\lesssim 
\sum_{\substack{\mu+|\alpha|\le M+1 \\ \mu\le \mu_0} }
\|L^\mu Z^\alpha u(0,\cdot)\|_{L^2(\ext)}
+
\sum_{\substack{\mu+|\alpha|\le M \\ \mu\le \mu_0} }
\|L^\mu \partial^\alpha u'\|_{L^2((0,T)\times\{x\in \ext : |x|<2\} ) },
\end{multline}
where we have used the Poincar{\'e} inequality for the last term.
Combining these estimates, we obtain the required result.
\ebox

\vspace{10pt}

Now, Lemma \ref{Lemma-WeightedEnergy-Ext} follows from Claim 
\ref{Lemma-EnergyKSSTangential-Claim2}, 
(2) of Lemma \ref{Lemma-Boundary-t}  
and the estimate
\begin{multline}
\sum_{\substack{\mu+|\alpha|\le M+2 \\ \mu\le \mu_0} }
\|L^\mu \partial^\alpha v\|_{L^2((0,T)\times\{x\in \whole : |x|<3\} ) }
\\
\lesssim 
\sum_{\substack{\mu+|\alpha|\le M+2 \\ \mu\le \mu_0} }
\|L^\mu \partial^\alpha u(0,\cdot)\|_{L^2(\ext)}
+
\sum_{\substack{\mu+|\alpha|\le M+1 \\ \mu\le \mu_0} }
\|L^\mu \partial^\alpha \square u\|_{L^1((0,T),L^2(\ext))},
\end{multline}
 which has been shown by Keel, Smith and Sogge \cite[(2.2), (2.5)]{Keel-Smith-Sogge-2002-JAnalMath}.
\ebox

%

\subsection{Sobolev type estimates}

We use the following weighted Sobolev estimates from 
\cite[Lemma 3.3]{Sogge-2008-InterPress}.  
To prove the estimate, we apply Sobolev estimates on $(0,\infty)\times S^2$.  
The decay result is from comparing the volume elements of $(0,\infty)\times S^2$ and $\mathbb{R}^3$.

\begin{lemma}
\label{Lemma-Sobolev-x}
Let $R\ge1$. The following inequality holds for any smooth function $h$. 
\begin{equation}
\|h\|_{L^\infty(\{x\in \whole : R/2<|x|<R\})} 
\lesssim 
R^{-1} \sum_{|\alpha|\le 2} \|Z^\alpha h\|_{L^2(\{ x\in \whole : R/4<|x|<2R\})}.
\end{equation}
\end{lemma}

\begin{lemma}
\label{Lemma-KS-Ext}
For any $M\ge0$ and $\mu_0\ge0$, the solution $u$ of \eqref{Eqn-Cauchy-F} satisfies the following.
\begin{multline}
\label{Ineq-KS-Ext}
(1) 
\sum_{\substack{\mu+|\alpha|\le M \\ \mu\le \mu_0}}
\angler^{1/2} \angletpr^{1/2} \anglectmr^{1/2}
|L^\mu Z^\alpha u'(t,x)|
\\
\lesssim
\sum_{\substack{\mu+|\alpha|\le M+2 \\ \mu\le \mu_0+1}}
\|L^\mu Z^\alpha u'(t,\cdot)\|_{L^2(\ext)}
+
\sum_{\substack{\mu+|\alpha|\le M+1 \\ \mu\le \mu_0}}
\|\langle t+|\cdot|\rangle L^\mu Z^\alpha \square_c u(t, \cdot)\|_{L^2(\ext )}
\\
+
\sum_{\substack{\mu\le \mu_0}}
t\|L^\mu u'(t,\cdot)\|_{L^2(\{y\in \ext : |y|<2\} )}.
\end{multline}

\begin{multline}
\label{Ineq-KS-Ext-v}
(2) 
\sum_{\substack{\mu+|\alpha|\le M \\ \mu\le \mu_0}}
\angler^{1/2} \angletpr^{1/2} \anglectmr^{1/2}
|L^\mu Z^\alpha u'(t,x)|
\\
\lesssim
\sum_{\substack{\mu+|\alpha|\le \mu_0+2 \\ \mu\le \mu_0}}
\|L^\mu \partial^\alpha u(0,\cdot)\|_{L^2(\ext)}
+
\sum_{\substack{\mu+|\alpha|\le M+2 \\ \mu\le \mu_0+1}}
\|L^\mu Z^\alpha u'(t,\cdot)\|_{L^2(\ext)}
\\
+
\sum_{\substack{\mu+|\alpha|\le \max\{M+1, \mu_0+1\} \\ \mu\le \mu_0}}
\sup_{0\le s\le t}
\|\langle s+|\cdot|\rangle L^\mu Z^\alpha \square_c u(s, \cdot)\|_{L^2(\ext )}
\\
+
\sum_{\substack{\mu+|\alpha|\le \mu_0+2 \\ \mu \le \mu_0}}
\sup_{0\le s\le t} (1+s)\|L^\mu \partial^\alpha v(t,\cdot)\|_{L^2( \{y\in \whole : |y|<3\} )}.
\end{multline}
\end{lemma}

\proof 
The proof of  (1) follows from  
\cite[Lemma 4.2, Lemma 4.3]{Metcalfe-Nakamura-Sogge-2005-ForumMath} 
and Lemma \ref{Lemma-KS}.
The proof of (2) follows from (1), and (1) of Lemma \ref{Lemma-Boundary-t}.


\subsection{Commutator estimates}
\label{Section-Commutator}

Let $\chi$ be a nonnegative function with $\chi(r)=0$ if $r\le 1$ and $\chi(r)=1$ if $r\ge2$.
We put 
\begin{equation}
\tilde{L}:=t\partial_t+\chi(r)r\partial_r.
\end{equation}

When we consider the version of higher derivatives by $L$ and $Z$ of  
\eqref{Eqn-Energy-Quasi}, we use the following commuting properties
(see 
\cite[p85, p113]{Metcalfe-Sogge-2005-InventMath}
and 
\cite[p4761]{Nakamura-2012-JDE}).

\begin{multline}
\label{Commutator-TildeL}
(\partial_t^2-c_I^2\Delta)\till^\mu\partial_t^m u_I
-
\sum_{\substack{1\le K\le D \\ 0\le k,l\le 3} } 
\gamma_{I}^{Kkl}\partial_k \partial_l \till^\mu \partial_t^m u_K
\\
=(\till+2)^\mu \partial_t^m \square_{\gamma_I} u_I
+
\sum_{\substack{p\le \mu-1 \\ |\nu|\le1} } \chi_{p\nu} \till^{p}\partial_t^m\partial^\nu \partial_x u_I
+
\sum_{K=1}^D
\sum_{\substack{p\le \mu-1 \\ |\nu|\le1 \\ 0\le k,l\le 3} }  
\chi_{p\nu klK} 
\gamma_{I}^{Kkl} \till^{p}\partial_t^m\partial^\nu \partial_x u_K
\\
+
\sum_{K=1}^D
\sum_{\substack{p+q\le \mu \\ r+s=m \\ p+r\ge1 \\ 0\le k,l\le 3} } 
C_{pqrsklK}\cdot 
(\till^p\partial_t^r \gamma_{I}^{Kkl}) \till^q\partial_t^s\partial_k\partial_l u_K,
\end{multline}
and 
\begin{multline}
\label{Commutator-L}
(\partial_t^2-c_I^2\Delta) L^\mu Z^\nu u_I
-
\sum_{\substack{1\le K\le D \\ 0\le k,l\le 3} } \gamma_{I}^{Kkl}\partial_k \partial_l L^\mu Z^\nu u_K
\\
=(L+2)^\mu Z^\nu \square_{\gamma_I} u_I 
+
\sum_{K=1}^D
\sum_{ \substack{p+q\le \mu \\ |\kappa|+|\delta|\le |\nu| \\ q+|\delta|\le \mu+|\nu|-1 \\ |\eta|=2} }
C_{pq\kappa\delta\eta K} \cdot(L^p Z^\kappa \gamma_{I}^{Kkl}) L^q Z^\delta \partial^\eta u_K,
\end{multline}
where $\chi_{p\nu}$ and $\chi_{p\nu klK}$ are smooth functions dependent on lower indices
 which supports are in the region $\{x\in \ext  : \chi(x)\le 2\}$, 
and the constants $C$ are dependent on the lower indices.

%

When we construct the energy estimates for the derivatives of the solution, we use the following estimates which follows from the elliptic regularity.
For any $M\ge0$ and $\mu_0\ge0$, we have 
\begin{multline}
\label{LtoTilL}
\sum_{\substack{\mu+|\alpha|\le M \\ \mu \le \mu_0} }
\|L^\mu\partial^\alpha u'(t,\cdot)\|_{L^2(\ext)}
\lesssim
\sum_{\substack{\mu+j\le M \\ \mu \le \mu_0} }
\|(\till^\mu\partial_t^j u)'(t,\cdot)\|_{L^2(\ext)}
\\
+
\sum_{\substack{\mu+|\alpha|\le M-1 \\ \mu \le \mu_0} }
\|L^\mu\partial^\alpha 
\square_c u(t,\cdot)\|_{L^2(\ext)}
\end{multline}
for any function $u$ and $c>0$ which satisfies the Dirichlet condition $u|_{\partial \mathcal{K}}=0$.

%
\subsection{Pointwise estimates}

We use the following pointwise estimates to show (4) of Proposition \ref{Prop-Ext}, below.
The first is the $L^\infty-L^\infty$ estimate due to Kubota and Yokoyama \cite[Theorem 3.4]{Kubota-Yokoyama-2001-JJM} for the boundaryless case (see also \cite[Theorem 2.4]{Metcalfe-Nakamura-Sogge-2005-JJM} and its proof).
The second is the $L^\infty-L^1$ estimate based on 
\cite[Proposition 2.1]{Keel-Smith-Sogge-2004-JAMS}.

\begin{lemma}
\label{Lemma-KY}
For any $\theta>0$, the solution $u=(u_1,\cdots,u_D)$ of 
\begin{equation}
\label{Eqn-Cauchy-F-Boundaryless}
\begin{array}{ll}
(\partial_t^2-c_I^2\Delta)u_I=F_I &\mbox{for}\ \ (t,x)\in [0,\infty)\times \mathbb{R}^3,\ \ 1\le I\le D \\
\end{array}
\end{equation}
satisfies 
\begin{multline}
(1+t+r)\left(1+\log\frac{1+t+r}{1+|c_It-r|}\right)^{-1}
|u_I(t,x)|
\\
\lesssim
\sum_{\substack{\mu+|\alpha|\le 3 \\ \mu\le 1, j\le 1} }
\|(\angley\partial)^jL^\mu Z^\alpha u_I(0,y)\|_{L^2_y}
\\
+
\sup_{(s,y)\in D_I(t,x)}|y|(1+s+|y|)^{1+\theta}\lambda_\theta(s,y) 
|F_I(s,y)|,
\end{multline}
where $r=|x|$ and $D_I$, $\lambda_\theta$ are defined by 
\begin{multline}
D_I(t,x)=\{(s,y)\in [0,\infty)\times\mathbb{R}^3 :
\\
0\le s\le t,
||x|-c_I(t-s)|\le |y|\le |x|+c_I(t-s)\},
\end{multline}
\begin{multline}
\Lambda_I=\{(s,y)\in [0,\infty)\times\mathbb{R}^3 :
\\
s\ge1,\ |y|\ge1,\ 
||y|-c_Is|\le \min_{1\le J,K\le D}|c_J-c_K|s/3\},
\end{multline}
\begin{multline}
\lambda_\theta(s,y)=
\left\{
\begin{array}{ll}
(1+||y|-c_Js|)^{1-\theta}  & \mbox{if}\ \ (s,y)\in \Lambda_J\ \ \mbox{for}\ \ 1\le J\le D \\
(1+|y|)^{1-\theta} & \mbox{if}\ \ (s,y)\in ((0,\infty)\times \mathbb{R}^3)\backslash (\cup_{1\le J\le D} \Lambda_J).
\end{array}
\right.
\end{multline}
\end{lemma}

\begin{lemma}
\label{Lemma-LInftyL1}
(\cite[Lemma 2.12]{Nakamura-2012-JDE})
Let $u$ be the solution of 
\begin{equation}
\label{Eq-Lemma-Huygens-Sublemma-1}
\left\{
\begin{array}{ll}
\square_c u=F &\mbox{for}\ \ (t,x)\in [0,\infty)\times \mathbb{R}^3 \\
u(t,x)=0 &\mbox{for}\ \ t\le0,\ \ x\in \mathbb{R}^3,
\end{array}
\right.
\end{equation}
where 
\[
\supp F\subset\{(s,y) : cs/10 \le |y| \le 10cs,\ \ s\ge1\}.
\]
Then 
\begin{multline}
(1+t+|x|)|u(t,x)|
\lesssim
\sum_{\substack{\mu+|\alpha|\le 3 \\ \mu\le 1} } 
\sup_{0\le s\le t} 
\int_{\mathbb{R}^3} |L^\mu Z^\alpha F(s,y)|dy
\left(
1+\left|\log\frac{1+t}{1+|c t-|x||}\right|
\right).
\end{multline}
\end{lemma}

%

\subsection{Estimates for nonlinear terms}
\label{Section-Nonlinear}

We show some estimates to treat the nonlinear terms.
We assume 
\begin{equation}
\label{Assume-NonlinearEstimates}
\sum_{|\alpha|\le M_0} \sup_{\substack{ 0\le t< T \\ x\in \ext} }
(1+t)|Z^\alpha u'(t,x)|\le C\varepsilon 
\end{equation}
for some constant $C>0$.
First we consider the semilinear part of quadratic nonlinearities. 
For $M\le 2M_0$, since we have 
\begin{multline}
\label{Ineq-Nonlinear-Modulo}
\sum_{\substack{\mu+|\alpha|\le M \\ \mu\le \mu_0 }} 
|L^\mu\partial^\alpha (u'u')|
\lesssim
\sum_{\substack{|\alpha|\le M_0 }} 
|\partial^\alpha u'|
\sum_{\substack{\mu+|\alpha|\le M \\ \mu\le \mu_0 }} 
|L^\mu\partial^\alpha u'|
\\
+
\sum_{\substack{M_0+1\le |\alpha|\le M-1 }} 
|\partial^\alpha u'|
\sum_{\substack{\mu+|\alpha|\le M-M_0-1 \\ 1\le \mu\le \mu_0 }} 
|L^\mu\partial^\alpha u'|
\\
+
\sum_{\substack{\mu+|\alpha|\le M/2 \\ 1\le \mu\le \mu_0-1}} 
|L^\mu\partial^\alpha u'|
\sum_{\substack{\mu+|\alpha|\le M-1 \\ 1\le \mu\le \mu_0-1 }} 
|L^\mu\partial^\alpha u'|.
\end{multline}
We obtain by the Sobolev estimate Lemma \ref{Lemma-Sobolev-x} and the assumption \eqref{Assume-NonlinearEstimates}
\begin{multline}
\label{Ineq-Nonlinear-Semi}
\sum_{\substack{\mu+|\alpha|\le M \\ \mu\le \mu_0 }} 
\|L^\mu\partial^\alpha (u'u')\|_{L^2(\ext)}
\lesssim
\frac{\varepsilon}{1+t}
\sum_{\substack{\mu+|\alpha|\le M \\ \mu\le \mu_0 }} 
\|L^\mu\partial^\alpha u'\|_{L^2(\ext)}
\\
+
\sum_{\substack{M_0+1\le |\alpha|\le M-1 }} 
\|\anglex^{-1/2} \partial^\alpha u'\|_{L^2(\ext)}
\sum_{\substack{\mu+|\alpha|\le M-M_0+1 \\ 1\le \mu\le \mu_0 }} 
\|\anglex^{-1/2} L^\mu Z^\alpha u'\|_{L^2(\ext)}
\\
+
\sum_{\substack{\mu+|\alpha|\le M/2+2 \\ 1\le \mu\le \mu_0-1}} 
\|\anglex^{-1/2} L^\mu Z^\alpha u'\|_{L^2(\ext)}
\sum_{\substack{\mu+|\alpha|\le M-1 \\ 1\le \mu\le \mu_0-1 }} 
\|\anglex^{-1/2} L^\mu\partial^\alpha u'\|_{L^2(\ext)}.
\end{multline}
Similarly, for the quasilinear part of the nonlinearities, we obtain for $M\le 2M_0-2$ 
\begin{multline}
\label{Ineq-Nonlinear-Quasi}
\sum_{\substack{\mu +|\alpha|+\nu+|\beta|\le M \\ \mu+\nu\le \mu_0 \\ \nu+|\beta|\le M-1 }} 
\|(L^\mu\partial^\alpha u') (L^\nu\partial^\beta u'')\|_{L^2(\ext)}
\lesssim
\frac{\varepsilon}{1+t}
\sum_{\substack{\mu+|\alpha|\le M \\ \mu\le \mu_0 }} 
\|L^\mu\partial^\alpha u'\|_{L^2(\ext)}
\\
+
\sum_{\substack{M_0+1\le |\alpha|\le M }} 
\|\anglex^{-1/2} \partial^\alpha u'\|_{L^2(\ext)}
\sum_{\substack{\mu+|\alpha|\le M-M_0+2 \\ 1\le \mu\le \mu_0 }} 
\|\anglex^{-1/2} L^\mu Z^\alpha u'\|_{L^2(\ext)}
\\
+
\sum_{\substack{\mu+|\alpha|\le M/2+3 \\ 1\le\mu\le \mu_0-1}} 
\|\anglex^{-1/2} L^\mu Z^\alpha u'\|_{L^2(\ext)}
\sum_{\substack{\mu+|\alpha|\le M \\ 1\le \mu\le \mu_0-1 }} 
\|\anglex^{-1/2} L^\mu\partial^\alpha u'\|_{L^2(\ext)}.
\end{multline}


\newsection{Continuity argument to prove Theorem \ref{Theorem-Ext} }
We prove the following proposition to prove Theorem \ref{Theorem-Ext}.

\begin{proposition}
\label{Prop-Ext}
Let $M_0\ge9$.
Let $\mathcal{K}$, $F$, $f$ and $g$ satisfy the assumption in Theorem \ref{Theorem-Ext}.
We put 
\begin{equation}
\varepsilon:=
\sum_{|\alpha|\le 2M_0}
\|\anglex^{|\alpha|+1} \partial_x^{\alpha} f(x)\|_{L^2(\ext)}
+
\sum_{|\alpha|\le 2M_0-1}
\|\anglex^{|\alpha|+1} \partial_x^{\alpha} g(x)\|_{L^2(\ext)}.
\end{equation}
Let $A_0>0$.
If the time local solution $u$ of the Cauchy problem \eqref{P} 
with the time interval $[0,\infty)$ replaced by $[0,T)$ satisfies 
\begin{equation}
\label{Prop-A0}
\sum_{\substack{|\alpha|\le M_0 \\ 1\le I\le D} }
\sup_{\substack{(t,x)\in S_T} }
\angler^{1/2}\angletpr^{1/2}\anglecItmr^{1/2}
|Z^\alpha u_I'(t,x)|  \le  A_0\varepsilon
\end{equation}
and $\varepsilon$ is sufficiently small, then for any $\mu_0\ge0$, $0\le M\le 2M_0-2-10\mu_0$ 
and $\sigma>0$,  
there exist constants $C_{M, \mu_0}>0$ which are dependent on $A_0$ such that the following inequality holds. 

\begin{multline}
\label{Prop-Higher}
\sum_{\substack{\mu+j\le M \\ \mu\le \mu_0}}
\|(\till^\mu\partial_t^j u)'(t,\cdot)\|_{L^2(\ext)}
+
\sum_{\substack{\mu+|\alpha|\le M \\ \mu\le \mu_0}}
\|L^\mu\partial^\alpha u'(t,\cdot)\|_{L^2(\ext)}
\\
+
\sum_{\substack{\mu+|\alpha|\le M-2 \\ \mu\le \mu_0}}
(\log(e+t))^{-1/2}
\|\langle x\rangle^{-1/2}L^\mu\partial^\alpha u'(\cdot,x)\|_{L^2(S_t)}
\\
+
\sum_{\substack{\mu+|\alpha|\le M-2 \\ \mu\le \mu_0 \\ 1\le I\le D}}
(\log(e+t))^{-1/2}
\|\langle c_Is-r \rangle^{-1/2}
\partialbar_{c_I} L^\mu\partial^\alpha u_I(s,x)\|_{L^2(S_t)}
\\
+
\sum_{\substack{\mu+|\alpha|\le M-3 \\ \mu\le \mu_0}}
\|L^\mu Z^\alpha u'(t,\cdot)\|_{L^2(\ext)}
\\
+
\sum_{\substack{\mu+|\alpha|\le M-5 \\ \mu\le \mu_0}}
(\log(e+t))^{-1/2}
\|\langle x\rangle^{-1/2}L^\mu Z^\alpha u'(\cdot,x)\|_{L^2(S_t)}
\\
+
\sum_{\substack{\mu+|\alpha|\le M-5 \\ \mu\le \mu_0 \\ 1\le I\le D}}
(\log(e+t))^{-1/2}
\|\langle c_Is-r\rangle^{-1/2}
\partialbar_{c_I} L^\mu Z^\alpha u_I(s,x)\|_{L^2(S_t)}
\\
\le 
C_{M,\mu_0}\varepsilon(1+t)^{C_{M,\mu_0}(\varepsilon+\sigma)}
\end{multline}
for $0\le t< T$.
Moreover, if $M_0\ge 32$, then there exist constants $C_0>0$, 
which is independent of $A_0$,  
and $C>0$, which is dependent on $A_0$, 
such that the following inequalities hold. 
\begin{equation}
\label{Ineq-Prop-Energy-High}
(1)
\sum_{\substack{\mu+|\alpha|\le M_0+5 \\ \mu\le 2}}
\|L^\mu Z^\alpha u'(t,\cdot)\|_{L^2(\ext)}
\le C_0\varepsilon 
+
C\varepsilon^2(1+t)^{C(\varepsilon+\sigma)},
\end{equation}
\begin{multline}
\label{Ineq-Prop-KS-High}
(2) 
\sum_{ \substack{\mu+|\alpha|\le M_0+3 \\ \mu\le 1 \\ 1\le I\le D} }
\sup_{x\in \ext}
\angler^{1/2}\angletpr^{1/2}\anglecItmr^{1/2}
|L^\mu Z^\alpha u_I'(t,x)|
\\
\le C_0\vep
+C\varepsilon^2(1+t)^{C(\varepsilon+\sigma)}
\end{multline}
for $0\le t< T$.
\begin{equation}
\label{Ineq-Prop-Energy-Lower}
\ \hspace{-2cm} 
(3) 
\sum_{\substack{\mu+|\alpha|\le M_0+2 \\ \mu\le 1}}
\|L^\mu Z^\alpha u'\|_{L^\infty((0,T),L^2(\ext))} 
\le
C_0\varepsilon+C\varepsilon^{3/2}.
\end{equation}
\begin{equation}
\label{Ineq-Prop-v}
\ \hspace{-2cm} 
(4) 
\sum_{\substack{ |\alpha|\le 2 }}
\sup_{0\le t< T}
(1+t)\|\partial^\alpha v(t,\cdot)\|_{L^\infty(\{x\in \whole : |x|<3\})} 
\le
C_0\varepsilon+C\varepsilon^{2},
\end{equation}
where $v=(v_1,\cdots,v_D)$ is the solution of 
\begin{equation}
\label{Eqn-Cauchy-vI}
\left\{
\begin{array}{l}
(\partial_t^2-c_I^2\Delta)v_I=(1-\zeta) F_I \ \ \mbox{for}\ \ (t,x)\in [0,T)\times \whole,\ \ 1\le I\le D \\
v(0,\cdot)=((1-\zeta) f)(\cdot),\ \ \partial_tv(0,\cdot)=((1-\zeta) g)(\cdot)
\end{array}
\right.
\end{equation}
and $\zeta \in C_0^\infty(\whole)$ is a function which satisfies 
$0\le \zeta \le1$, 
$\zeta (x)=1$ for $|x|\le 3$, and $\zeta (x)=0$ for $|x|\ge4$.
Here, $(1-\zeta) F$, $(1-\zeta) f$, $(1-\zeta) g$ are regarded as functions on $\whole$ by zero-extension.

\begin{multline}
\label{Ineq-Sobolev-Lower}
(5) 
\sum_{\substack{|\alpha|\le M_0 \\ 1\le I\le D }}
\sup_{\substack{(t,x)\in S_T} }
\angler^{1/2} \angletpr^{1/2} \langle c_It-r\rangle^{1/2}
|Z^\alpha u_I'(t,x)|
\le C_0\varepsilon+C\varepsilon^{3/2}.
\end{multline}

\end{proposition}


\subsection{Proof of Theorem \ref{Theorem-Ext}}
We use the continuity argument which shows that the local in time solution $u$ does not blow up if its initial data are sufficiently small.
We refer to 
\cite{Keel-Smith-Sogge-2002-JFA} 
for the existence of the local in time solutions.
Since the constant $C_0$ is independent of $A_0$ in \eqref{Ineq-Sobolev-Lower}, 
we put $A_0=4C_0$ and take $\varepsilon$ sufficiently small such that 
$C\varepsilon^{3/2}\le C_0\varepsilon$.
Then the right hand side of \eqref{Ineq-Sobolev-Lower} is bounded by 
$A_0\varepsilon/2$, which shows the local in time solution $u$ does not blow up, 
namely the solution exists globally in time.
\ebox

%

\subsection{Proof of Proposition \ref{Prop-Ext}}

First, we show the estimate \eqref{Prop-Higher} inductively, and then we derive the estimates from (1) to (5).
We drop the indices $I$ of $c_I$, $u_I$ and so on to avoid the complexity.

\subsubsection{The estimate for $\|L^\mu \partial^\alpha u'\|_2$}

By \eqref{LtoTilL}, we have  
\begin{equation}
\sum_{\substack{\mu+|\alpha|\le M \\ \mu\le \mu_0 }} 
\|L^\mu\partial^\alpha u'\|_2
\lesssim
\sum_{\substack{\mu+j\le M \\ \mu\le \mu_0 }}
\|\partial \till^\mu\partial_t^ju\|_2
+
\sum_{\substack{\mu+|\alpha|\le M-1 \\ \mu\le \mu_0 }}
\|L^\mu\partial^\alpha \square_c u\|_2,
\end{equation}
and we estimate the last term by the argument in Section  
\ref{Section-Nonlinear}
\begin{equation}
\sum_{\substack{\mu+|\alpha|\le M-1 \\ \mu\le \mu_0 }} 
\|L^\mu\partial^\alpha \square_c u\|_2
\lesssim
\frac{\varepsilon}{1+t}
\sum_{\substack{\mu+|\alpha|\le M \\ \mu\le \mu_0 }} 
\|L^\mu\partial^\alpha u'\|_2
+A,
\end{equation}
where 
\begin{multline}
A:=
\sum_{\substack{M_0+1\le |\alpha|\le M-1 } } 
\|\anglex^{-1/2} \partial^\alpha u'\|_2
\sum_{\substack{\mu+|\alpha|\le M-M_0+1 \\ 1\le\mu\le \mu_0 }} 
\|\anglex^{-1/2} L^\mu Z^\alpha u'\|_2
\\
+
\sum_{\substack{\mu+|\alpha|\le M/2+2 \\ 1\le\mu\le \mu_0-1 }} 
\|\anglex^{-1/2} L^\mu Z^\alpha u'\|_2
\sum_{\substack{\mu+|\alpha|\le M-1 \\ 1\le\mu\le \mu_0-1 }} 
\|\anglex^{-1/2} L^\mu \partial^\alpha u'\|_2.
\end{multline}
Therefore for sufficiently small $\varepsilon>0$, we obtain 
\begin{equation}
\label{Ineq-Till-L}
\sum_{\substack{\mu+|\alpha|\le M \\ \mu\le \mu_0 }} 
\|L^\mu\partial^\alpha u'\|_2
\lesssim
\sum_{\substack{\mu+j\le M \\ \mu\le \mu_0 }}
\|\partial \till^\mu\partial_t^ju\|_2
+A.
\end{equation}

\subsubsection{The estimate for the boundary term}

To consider the estimate for $\|\partial \till^\mu\partial_t^j u\|_2$ in the next subsection, 
we prepare the estimate for the boundary term.
By Lemma \ref{Lemma-BoundaryTerm}, we have
\begin{multline}
\sum_{\substack{\mu+|\alpha|\le M \\ \mu\le \mu_0-1 }} 
\int_0^t\|L^\mu\partial^\alpha u'(s,x)\|_{L^2(|x|<2)}ds
\lesssim
\sum_{\substack{\mu+|\alpha|\le M+2 \\ \mu\le \mu_0-1 }} 
\|\anglex L^\mu \partial^\alpha u(0,x)\|_{L^2}
\\
+
\sum_{\substack{\mu+|\alpha|\le M+1 \\ \mu\le \mu_0-1 }} 
\int_0^t \|L^\mu\partial^\alpha \square_c u(s,x)\|_{L^2_x(|x|<4)}ds
\\
+
\sum_{\substack{\mu+|\alpha|\le M+1 \\ \mu\le \mu_0-1 }} 
\int_0^t \int_0^s \|L^\mu\partial^\alpha \square_c u(\tau,x)\|_{L^2_x(||x|-c(s-\tau)|<4)}d\tau ds.
\end{multline}
Since the last two terms are bounded by  
\begin{equation}
\sum_{\substack{\mu+|\alpha|\le \max\{M/2+3, M+2\} \\ \mu\le \mu_0-1}}
\|\anglex^{-1/2} L^\mu Z^\alpha u'(s,x)\|_{L^2_{s,x}(S_t)}^2
\end{equation}
due to  Lemma \ref{Lemma-Sobolev-x}, we obtain 
\begin{multline}
\label{Ineq-BoundaryTerms}
\sum_{\substack{\mu+|\alpha|\le M \\ \mu\le \mu_0-1 }} 
\int_0^t\|L^\mu\partial^\alpha u'\|_{L^2(|x|<2)}ds
\lesssim
\sum_{\substack{\mu+|\alpha|\le M+2 \\ \mu\le \mu_0-1 }} 
\|\anglex L^\mu \partial^\alpha u(0,x)\|_{L^2_x(\ext)}
\\
+
\sum_{\substack{\mu+|\alpha|\le \max\{M/2+3, M+2\} \\ \mu\le \mu_0-1}}
\|\anglex^{-1/2} L^\mu Z^\alpha u'\|_{L^2(S_t)}^2.
\end{multline}

%

\subsubsection{The estimate for $\|\partial \till^\mu \partial_t^ju\|_2$}

Since $\till^\mu \partial_t^ju$ satisfies the Dirichlet condition, by the energy estimate \eqref{Eqn-Divergence}, we have 
\begin{multline}
\partial_t 
\Big\{
\int_\ext e_0(\till^\mu \partial_t^ju)dx
\Big\}^{1/2}
\lesssim
\|\square_\gamma \till^\mu \partial_t^ju\|_2
+
\|\gamma'\|_\infty
\Big\{
\int_\ext e_0(\till^\mu \partial_t^ju)dx
\Big\}^{1/2}.
\end{multline}
By \eqref{Commutator-TildeL}, we have 
\begin{multline}
\sum_{\substack{\mu+j\le M \\ \mu\le \mu_0 }} 
\|\square_\gamma \tilde{L}^\mu\partial_t^j u\|_{L^2(\ext)}
\lesssim
(1+\|\gamma\|_\infty)
\sum_{\substack{\mu+|\alpha|\le M \\ \mu\le \mu_0-1 }} 
\|L^\mu \partial^\alpha u'\|_{L^2(|x|<2)}
\\
+
\sum_{\substack{\mu+|\alpha|\le M \\ \mu\le \mu_0 }} 
\|L^\mu \partial^\alpha \square_\gamma u\|_2
+
\sum_{\substack{\mu+|\alpha|+\nu+|\beta|\le M \\ \mu+\nu\le \mu_0 
\\ \nu+|\beta|\le M-1 }} 
\|(L^{\mu}\partial^{\alpha}\gamma)
\cdot
(L^{\nu}\partial^{\beta}u'')\|_2.
\end{multline}
Using the estimates \eqref{Ineq-Nonlinear-Semi} and \eqref{Ineq-Nonlinear-Quasi},
the last two terms are bounded by 
\begin{equation}
\frac{\varepsilon}{1+t}
\sum_{\substack{\mu+|\alpha|\le M \\ \mu\le \mu_0 }} 
\|L^\mu\partial^\alpha u'\|_2
+B 
\end{equation}
for $M\le 2M_0-2$, where 
\begin{multline}
\label{Def-B}
B:=
\sum_{\substack{M_0+1\le |\alpha|\le M  }} 
\|\langle x\rangle^{-1/2}\partial^\alpha u'\|_2
\sum_{\substack{\mu+|\alpha|\le M-M_0+2  \\ 1\le \mu\le \mu_0}} 
\|\langle x\rangle^{-1/2}L^\mu Z^\alpha u'\|_2
\\
+
\sum_{\substack{\mu+|\alpha|\le M/2+3  \\ 1\le \mu\le \mu_0-1}} 
\|\langle x\rangle^{-1/2}L^\mu Z^\alpha u'\|_2
\sum_{\substack{\mu+|\alpha|\le M  \\ 1\le \mu\le \mu_0-1}} 
\|\langle x\rangle^{-1/2}L^\mu\partial^\alpha u'\|_2.
\end{multline}
So that, by \eqref{Ineq-Till-L}, the Gronwall inequality and \eqref{Ineq-BoundaryTerms}, we have 
\begin{multline}
\sum_{\substack{\mu+j\le M \\ \mu\le \mu_0 }} 
\|\partial \till^\mu \partial_t^ju\|_{L^2(\ext)}
\lesssim 
\sum_{\substack{\mu+j\le M \\ \mu\le \mu_0 }} 
\Big\{
\int_\ext e_0(\till^\mu \partial_t^ju)dx
\Big\}^{1/2}
\\
\lesssim
\Big\{
\sum_{\substack{\mu+|\alpha|\le M+2 \\ \mu\le \mu_0 }} 
\|\anglex L^\mu\partial^\alpha u(0,x)\|_{L^2_x(\ext)}
+
\sum_{\substack{\mu+|\alpha|\le \max\{M/2+3,M+2\} \\ \mu\le \mu_0-1 }} 
\|\anglex^{-1/2} L^\mu Z^\alpha u'\|_{L^2(S_t)}^2
\\
+
\sum_{\substack{M_0+1\le |\alpha|\le M} } 
\|\anglex^{-1/2}\partial^\alpha u'\|_{L^2(S_t)}
\sum_{\substack{\mu+|\alpha|\le M-M_0+2 \\ 1\le \mu\le \mu_0 }} 
\|\anglex^{-1/2}L^\mu Z^\alpha u'\|_{L^2(S_t)}
\\
+
\sum_{\substack{\mu+|\alpha|\le M/2+3 \\ 1\le \mu\le \mu_0-1 }} 
\|\anglex^{-1/2}L^\mu Z^\alpha u'\|_{L^2(S_t)}
\sum_{\substack{\mu+|\alpha|\le M \\ 1\le \mu\le \mu_0-1 }} 
\|\anglex^{-1/2}L^\mu \partial^\alpha u'\|_{L^2(S_t)}
\Big\}(1+t)^{C\varepsilon}
\end{multline}
for $M\le 2M_0-2$, where we have used that $\varepsilon>0$ is sufficiently small for the first inequality.


\subsubsection{The estimate for $\|L^\mu Z^\alpha u'\|_2$}

By the energy estimate \eqref{Eqn-Divergence} for $L^\mu Z^\alpha u$, we have 
\begin{multline}
\label{Energy-Ext-LZu}
\sum_{\substack{\mu+|\alpha|\le M \\ \mu\le \mu_0 }} 
\partial_t\int_\ext e_0(L^\mu Z^\alpha u) dx
\lesssim
\sum_{\substack{\mu+|\alpha|\le M+1 \\ \mu\le \mu_0 }} 
\|L^\mu \partial^\alpha u'\|_{L^2(\{x\in\ext : |x|<2\})}^2
\\
+
\sum_{\substack{\mu+|\alpha|\le M \\ \mu\le \mu_0 }} 
\left|
\int_\ext (\partial_tL^\mu Z^\alpha u) \square_\gamma L^\mu Z^\alpha udx
\right|
+
\|\gamma'\|_\infty 
\sum_{\substack{\mu+|\alpha|\le M \\ \mu\le \mu_0 }} 
\|L^\mu Z^\alpha u'\|_2^2,
\end{multline}
where we have used the trace theorem for the boundary term, so that, there is a loss of one derivative.
By \eqref{Commutator-L}, we have 
\begin{multline}
\sum_{\substack{\mu+|\alpha|\le M \\ \mu\le \mu_0}}
|\square_\gamma L^\mu Z^\alpha u|
\lesssim
\sum_{\substack{\mu+|\alpha|\le M \\ \mu\le \mu_0 }} 
|L^\mu Z^\alpha \square_\gamma u|
\\
+
\sum_{\substack{\mu+|\alpha|+\nu+|\beta|\le M \\ \mu+\nu\le \mu_0 \\ \nu+|\beta|\le M-1}} 
|(L^\mu Z^\alpha \gamma)(L^\nu Z^\beta \partial^2 u)|.
\end{multline}
So that, by \eqref{Ineq-Nonlinear-Semi} and \eqref{Ineq-Nonlinear-Quasi}, 
we have 
\begin{equation}
\sum_{\substack{\mu+|\alpha|\le M \\ \mu\le \mu_0 }} 
\|\square_\gamma L^\mu Z^\alpha u\|_2
\lesssim
\frac{\varepsilon}{1+t}
\sum_{\substack{\mu+|\alpha|\le M \\ \mu\le \mu_0 }} 
\|L^\mu Z^\alpha u'\|_2+\tilde{B},
\end{equation}
where $\tilde{B}$ is $B$ in \eqref{Def-B} with all $\partial$ replaced by $Z$.
Since $\|L^\mu Z^\alpha u'\|_2$ is equivalent to $e_0(L^\mu Z^\alpha u)$ for sufficiently small $\varepsilon>0$, we have 
\begin{multline}
\sum_{\substack{\mu+|\alpha|\le M \\ \mu\le \mu_0 }} 
\partial_t\int_\ext e_0(L^\mu Z^\alpha u) dx
\lesssim
\sum_{\substack{\mu+|\alpha|\le M+1 \\ \mu\le \mu_0 }} 
\|L^\mu \partial^\alpha u'\|_{L^2(|x|<2)}^2
\\
+
\frac{\varepsilon}{1+t}
\sum_{\substack{\mu+|\alpha|\le M \\ \mu\le \mu_0 }} 
\int_\ext e_0(L^\mu Z^\alpha u)dx
+
\left(
\sum_{\substack{\mu+|\alpha|\le M \\ \mu\le \mu_0 }} 
\int_\ext e_0(L^\mu Z^\alpha u)dx
\right)^{1/2}\cdot\tilde{B}.
\end{multline}
So that, by the Gronwall inequality, we obtain for $M\le 2M_0-2$
\begin{multline}
\sum_{\substack{\mu+|\alpha|\le M \\ \mu\le \mu_0 }} 
\|L^\mu Z^\alpha u'\|_{L^2(\ext)}
\lesssim
\sum_{\substack{\mu+|\alpha|\le M \\ \mu\le \mu_0 }} 
\left\{
\int_\ext e_0(L^\mu Z^\alpha u) dx
\right\}^{1/2}
\\
\lesssim
\Big\{
\sum_{\substack{\mu+|\alpha|\le M \\ \mu\le \mu_0 }} 
\|L^\mu Z^\alpha u'(0,\cdot)\|_{L^2(\ext)}
+
\sum_{\substack{\mu+|\alpha|\le M+1 \\ \mu\le \mu_0 }} 
\|\anglex^{-1/2} L^\mu \partial^\alpha u'\|_{L^2(S_t)}
\\
+
\sum_{\substack{M_0+1\le |\alpha|\le M }} 
\|\anglex^{-1/2}Z^\alpha u'\|_{L^2(S_t)}
\sum_{\substack{\mu+ |\alpha|\le M-M_0+2 \\ 1\le \mu\le \mu_0 }} 
\|\anglex^{-1/2}L^\mu Z^\alpha u'\|_{L^2(S_t)}
\\
+
\sum_{\substack{\mu+|\alpha|\le M/2+3 \\ 1\le \mu\le \mu_0-1 }} 
\|\anglex^{-1/2}L^\mu Z^\alpha u'\|_{L^2(S_t)}
\sum_{\substack{\mu+ |\alpha|\le M \\ 1\le \mu\le \mu_0-1 }} 
\|\anglex^{-1/2}L^\mu Z^\alpha u'\|_{L^2(S_t)}
\Big\}
(1+t)^{C\varepsilon}.
\end{multline}

%

\subsubsection{The estimates for the weighted energy}

By Lemma \ref{Lemma-WeightedEnergy-Ext}  
and \eqref{Ineq-Nonlinear-Semi}, 
we have for $M\le 2M_0-2$
\begin{multline}
\label{Ineq-WeightedEnergy-High}
\sum_{\substack{\mu+|\alpha|\le M \\ \mu\le \mu_0 }} 
\|L^\mu \partial^\alpha u'\|_{L^\infty((0,t),L^2(\ext))}
+
\sum_{\substack{\mu+|\alpha|\le M \\ \mu\le \mu_0 }} 
(\log(e+t))^{-1/2}
\|\anglex^{-1/2}L^\mu \partial^\alpha u'\|_{L^2(S_t)}
\\
+
\sum_{\substack{\mu+|\alpha|\le M \\ \mu\le \mu_0 \\ 1\le I\le D}}
(\log(e+t))^{-1/2}
\|\langle c_Is-r \rangle^{-1/2}
\partialbar_{c_I} L^\mu\partial^\alpha u_I(s,x)\|_{L^2(S_t)}
\\
\lesssim
\sum_{\substack{\mu+|\alpha|\le M+2 \\ \mu\le \mu_0 }} 
\|L^\mu \partial^\alpha u(0,\cdot)\|_{L^2(\ext)}
+
\sum_{\substack{\mu+|\alpha|\le M+2 \\ \mu\le \mu_0 }} 
\|L^\mu \partial^\alpha u'\|_{L^\infty((0,t),L^2(\ext))}
\varepsilon\log(1+t)
\\
+
\sum_{\substack{M_0+1\le |\alpha|\le M+1 }} 
\|\anglex^{-1/2}\partial^\alpha u'\|_{L^2(S_t)}
\sum_{\substack{\mu+ |\alpha|\le M-M_0+3 \\ 1\le \mu\le \mu_0 }} 
\|\anglex^{-1/2}L^\mu Z^\alpha u'\|_{L^2(S_t)}
\\
+
\sum_{\substack{\mu+|\alpha|\le M/2+3 \\ 1\le \mu\le \mu_0-1 }} 
\|\anglex^{-1/2}L^\mu Z^\alpha u'\|_{L^2(S_t)}
\sum_{\substack{\mu+ |\alpha|\le M+1 \\ 1\le \mu\le \mu_0-1 }} 
\|\anglex^{-1/2}L^\mu \partial^\alpha u'\|_{L^2(S_t)}.
\end{multline}
Here, the above estimate also holds with all $\partial$ replaced by $Z$.

%

\subsubsection{The proof of (1)}
The proof of (1) follows from \eqref{Ineq-WeightedEnergy-High} 
with $M=M_0+5$ and $\mu_0=2$ by \eqref{Prop-Higher}.
Indeed, to bound 
$\sum_{\substack{\mu+|\alpha|\le M_0+5 \\ \mu\le 2 }} 
\|L^\mu Z^\alpha u'\|_{L^\infty((0,t),L^2(\ext))}$, 
we need the estimate for 
$\sum_{\substack{\mu+|\alpha|\le M_0+7 \\ \mu\le 2 }} 
\|L^\mu Z^\alpha u'\|_{L^\infty((0,t),L^2(\ext))}$ 
by \eqref{Ineq-WeightedEnergy-High},
which is bounded by $C\varepsilon(1+t)^{C(\varepsilon+\sigma)}$ by \eqref{Prop-Higher} since 
$M_0+7\le (2M_0-2-20)-3$ is satisfied by $M_0\ge32$.

%

\subsubsection{The proof of (2)}

By \eqref{Prop-A0} and induction argument, we have 
\begin{multline}
\sum_{\substack{\mu+|\alpha|\le M_0+4 \\ \mu\le 1}}
\langle t+|y| \rangle |L^\mu Z^\alpha \square_c u(t,y)|
\lesssim
\varepsilon 
\sum_{\substack{\mu+|\alpha|\le M_0+5 \\ \mu\le 1}}
|L^\mu Z^\alpha u'(t,y)|(1+t)^{C(\varepsilon+\sigma)}.
\end{multline}
So that, we have  
\begin{multline}
\sum_{\substack{\mu+|\alpha|\le M_0+4 \\ \mu\le 1}}
\|\langle t+|y| \rangle L^\mu Z^\alpha \square_c u(t,y)\|_{L^2_y(\ext)}
\\
\lesssim
\varepsilon 
\sum_{\substack{\mu+|\alpha|\le M_0+5 \\ \mu\le 1}}
\|L^\mu Z^\alpha u'(t,y)\|_{L^2_y(\ext)}(1+t)^{C(\varepsilon+\sigma)}.
\end{multline}
By simple calculation and \eqref{Ineq-BoundaryTerms}, we have 
\begin{multline}
\sum_{\substack{\mu \le 1}}
t\|L^\mu u'(t,y)\|_{L^2(|y|<2)}
\lesssim
\sum_{\substack{\mu+|\alpha|\le 2 \\ \mu\le 1}}
\int_0^t \|L^\mu\partial^\alpha u'(s,y)\|_{L^2(|y|<2)} ds
\\
\lesssim 
\sum_{\substack{\mu+|\alpha|\le 4 \\ \mu\le 1}}
\|\angley L^\mu\partial^\alpha u(0,y)\|_{L^2_y(\ext)} 
+
\sum_{\substack{\mu+|\alpha|\le 4 \\ \mu\le 1}}
\|\anglex^{-1/2} L^\mu Z^\alpha u'\|_{L^2(S_t)}^2.
\end{multline}
So that, by \eqref{Ineq-KS-Ext}, we obtain 
\begin{multline}
\sum_{\substack{\mu+|\alpha|\le M_0+3 \\ \mu\le 1}}
\angler^{1/2} \angletpr^{1/2} \anglectmr^{1/2}
|L^\mu Z^\alpha u'(t,x)|
\\
\lesssim
\sum_{\substack{\mu+|\alpha|\le 4 \\ \mu\le 1}}
\|\angley L^\mu\partial^\alpha u(0,y)\|_{L^2_y(\ext)} 
+
\sum_{\substack{\mu+|\alpha|\le 4 \\ \mu\le 1}}
\|\anglex^{-1/2} L^\mu Z^\alpha u'\|_{L^2(S_t)}^2
\\
+
\varepsilon 
\sum_{\substack{\mu+|\alpha|\le M_0+5 \\ \mu\le 1}}
\|L^\mu Z^\alpha u'(t,y)\|_{L^2_y(\ext)}(1+t)^{C(\varepsilon+\sigma)}
\\
+
\sum_{\substack{\mu+|\alpha|\le M_0+5 \\ \mu\le 2}}
\|L^\mu Z^\alpha u'(t,\cdot)\|_{L^2(\ext)}.
\end{multline}
Therefore we obtain the required estimate by \eqref{Prop-Higher} 
and (1).

%

\subsubsection{The proof of (3)}

By the standard energy estimate, we have 
\begin{multline}
\sum_{\substack{\mu+|\alpha|\le M_0+2 \\ \mu\le 1}}
\|L^\mu Z^\alpha u_I'(t,\cdot)\|_{L^2}^2
\le
C_0 \sum_{\substack{\mu+|\alpha|\le M_0+2 \\ \mu\le 1}}
\|L^\mu Z^\alpha u_I'(0,\cdot)\|_{L^2}^2 
\\
+
C_0 \sum_{\substack{\mu+|\alpha|\le M_0+2 \\ \mu\le 1}}
\int_0^t\int_{\ext} 
|(\partial_t L^\mu Z^\alpha u_I )
\square_{c_I} L^\mu Z^\alpha u_I| dxds
\\
+
C_0 \sum_{\substack{\mu+|\alpha|\le M_0+3 \\ \mu\le 1}}
\int_0^t \|L^\mu \partial^\alpha u_I'(s,x)\|_{L^2(|x|<1)}^2ds, 
\end{multline}
where $C_0>0$ is independent of $A_0$.
We use (2) to bound the last term 
\begin{multline}
\sum_{\substack{\mu+|\alpha|\le M_0+3 \\ \mu\le 1}}
\int_0^t \|L^\mu Z^\alpha u_I'(s,x)\|_{L^2(|x|<1)}^2ds
\\
\le
\int_0^t \|\angler^{-1/2}\langle s+r\rangle^{-1/2}\langle c_Is-r\rangle^{-1/2}
\{C_0\varepsilon+C\varepsilon^2(1+s)^{C(\varepsilon+\sigma)}\}\|_{L^2(|x|<1)}^2 ds
\\
\le
(C_0\varepsilon+C\varepsilon^2)^2.
\end{multline}
Since we are able to have the bound 
\beq
\sum_{\substack{\mu+|\alpha|\le M_0+2 \\ \mu\le 1}}
\int_0^t\int_{\ext} 
|(\partial_t L^\mu Z^\alpha u_I )
\square_{c_I} L^\mu Z^\alpha u_I| dxds
\le C\varepsilon^3
\eeq
by the similar argument for the proof of \eqref{Ineq-LowEnergy}, we obtain the required inequality.
%
\subsubsection{The proof of (4)}

Let $\chi\in C^\infty(\br)$ satisfy $\chi(t)=0$ for $t\le 1$, 
and $\chi(t)=1$ for $t\ge2$.
Let $\eta\in C^\infty(\br)$ satisfy $\eta(r)=0$ for 
$r\le \min_{1\le I\le D} c_I/10$ or 
$r\ge 10\max_{1\le I\le D} c_I$, and  $\eta(r)=1$ for 
$\min_{1\le I\le D} c_I/5\le r\le 5\max_{1\le I\le D} c_I$.
We put $\rho(t,x):=\chi(t)\eta(|x|/t)$.
We decompose $v$ into $w=(w_1,\cdots,w_D)$ and $z=(z_1,\cdots,z_D)$ which satisfy
\beq
\left\{
\begin{array}{l}
\square_{c_I} w_I(t,x)=\rho(t,x)\square_{c_I} v_I(t,x)\ \ \mbox{for}\ \ (t,x)\in [0,T)\times\whole,\ \ 1\le I\le D \\
w(0,\cdot)=\partial_t w(0,\cdot)=0
\end{array}
\right.
\eeq
\beq
\left\{
\begin{array}{l}
\square_{c_I} z_I(t,x)=(1-\rho(t,x))\square_{c_I} v_I(t,x)\ \ \mbox{for}\ \ (t,x)\in [0,T)\times\whole,\ \ 1\le I\le D \\
z(0,\cdot)=v(0,\cdot),\ \ \partial_t z(0,\cdot)=\partial_tv (0,\cdot).
\end{array}
\right.
\eeq
We note $v=w+z$ and show the required estimates for $w$ and $z$.
By Lemma \ref{Lemma-LInftyL1}, we have
\begin{multline}
\sum_{\substack{|\alpha|\le 2} }
(1+t+|x|)|Z^\alpha w_I(t,x)|
\\
\lesssim
\sum_{\substack{\mu+|\alpha|\le 5 \\ \mu\le 1} }
\sup_{0\le s\le t}\int_{\whole} |L^\mu Z^\alpha \square_{c_I} w_I(s,y)|dy
\left(
1+\left|\log\frac{1+t}{1+|c_It-|x|} \right|
\right).
\end{multline}
So that, we have by (3) 
\beq
\sum_{\substack{|\alpha|\le 2} }
(1+t)\|\partial^\alpha w_I(t,x)\|_{L^\infty(|x|<3)}
\lesssim
\sum_{\substack{\mu+|\alpha|\le 6 \\ \mu\le 1} }
\|L^\mu Z^\alpha u'\|_{L^\infty((0,t), L^2(\ext))}^2
\lesssim \varepsilon^2.
\eeq

On the other hand, by Lemma \ref{Lemma-KY}, we have 
\begin{multline}
\sum_{\substack{|\alpha|\le 2 } }
(1+t+|x|)|Z^\alpha z_I(t,x)|
\lesssim
\sum_{\substack{\mu+|\alpha|\le 5 \\ \mu\le 1, j\le 1} }
\|(\angley\partial)^j L^\mu Z^\alpha z_I(0,y)\|_{L^2_y(\whole)}
\\
+
\sum_{\substack{|\alpha|\le 2 } }
\sup_{(s,y)\in D_I(t,x)} |y| (1+s+|y|)^{1+\theta} \lambda_\theta(s,y) 
|\square_{c_I} Z^\alpha z_I(s,y)|
\end{multline}
for any fixed $\theta>0$. Since
\begin{multline}
\sum_{\substack{|\alpha|\le 2 } }
|\square_{c_I} Z^\alpha z_I(s,y)|
\lesssim
(1-\rho(s,y))
\sum_{\substack{|\alpha|\le 3} }
|Z^\alpha u'(s,y)|^2
\lesssim
\varepsilon^2\angley^{-1}\langle s+|y| \rangle^{-2}
\end{multline}
by \eqref{Prop-A0}, we obtain 
\beq
\sum_{\substack{|\alpha|\le 2 } }
(1+t+|x|)|Z^\alpha z_I(t,x)|
\le
C_0\varepsilon +C\varepsilon^2.
\eeq
Combining the above estimates for $w_I$ and $z_I$, we obtain the required estimate.

\subsubsection{The proof of (5)}

By \eqref{Ineq-KS-Ext-v}, we have 
\begin{multline}
\sum_{\substack{|\alpha|\le M_0}}
\angler^{1/2} \angletpr^{1/2} \anglectmr^{1/2}
|Z^\alpha u'(t,x)|
\le
C_0\sum_{\substack{|\alpha|\le 2}}
\|\partial^\alpha u(0,\cdot)\|_{L^2(\ext)}
\\
+
C_0\sum_{\substack{\mu+|\alpha|\le M_0+2 \\ \mu\le 1}}
\|L^\mu Z^\alpha u'(t,\cdot)\|_{L^2(\ext)}
+
C_0\sum_{\substack{|\alpha|\le M_0+1 }}
\sup_{0\le s\le t}
\|\langle s+|\cdot|\rangle Z^\alpha \square_c u(s, \cdot)\|_{L^2(\ext )}
\\
+
C_0\sum_{\substack{|\alpha|\le 2}}
\sup_{0\le s\le t} (1+s)
\|\partial^\alpha v(t,\cdot)\|_{L^2( \{y\in \whole : |y|<3\} )}
=: D_1+D_2+D_3+D_4.
\end{multline}
Since we have 
\begin{equation}
D_3\lesssim
\varepsilon \sum_{|\alpha|\le M_0+2} \|Z^\alpha u'\|_{L^\infty((0,t),L^2)}
\end{equation}
by \eqref{Prop-A0}, we obtain the required result by (3) and (4).
\ebox

%

\newsection{Appendices}
We put two notes on the weighted energy estimates and the wave equations with single speed.

\subsection{Weighted energy estimates}
We prove the following lemma which generalizes the weighted energy estimate of tangential derivatives in Lemma 2.1.
Let $n\ge1$, $c>0$, $\Delta:=\sum_{j=1}^n \partial_j^2$ 
and $\nabla:=(\partial_1,\cdots,\partial_n)$.
We denote the tangential derivatives along the $c$ speed light cone by  
\beq
\label{Def-Partialbar}
\partialbar_c
=(\partialbar_{c0},\partialbar_{c1}, \cdots,\partialbar_{cn})
:=
\left\{
\begin{array}{l}
(\partial_t+c\partial_r, \nabla-\frac{x}{r}\partial_r)\ \ \mbox{for}\ \ n\ge2 \\
(\partial_t+c\partial_r,0)\ \ \mbox{for}\ \ n=1,
\end{array}
\right.
\eeq
where $r:=|x|$ and $\partial_r:=r^{-1}(\sum_{1\le j\le n}x_j\partial_j)$.
For any $f$, $g$ and $F$, we consider the Cauchy problem
\beq
\label{Cauchy-n}
\left\{
\begin{array}{l}
(\partial_t^2-c^2\Delta)u(t,x)=F(t,x)\ \ \ \ \mbox{for}\ \ (t,x)\in [0,T]\times\mathbb{R}^n \\
u(0,\cdot)=f(\cdot), \ \ \partial_tu(0,\cdot)=g(\cdot).
\end{array}
\right.
\eeq 

\begin{lemma}
\label{Lemma-WE}
Let $n\ge1$.
The solution $u$ of \eqref{Cauchy-n} satisfies the following estimate.
\begin{multline}
\label{WeightedEnergy-Tangential-n}
\max\Big[
\sup_{0\le t\le T}
\int_{\mathbb{R}^n}
(\partial_t u(t,x))^2+c^2
|\nabla u(t,x)|^2
dx,
\\
\sup_{0<\kappa<\infty}
\frac{c\kappa}{4}
\int_0^T \int_{\mathbb{R}^n}
\frac{1}{ (1+|ct-|x||)^{1+\kappa} } 
\left\{
(\partialbar_{c0}u)^2+c^2\sum_{j=1}^n (\partialbar_{cj}u)^2
\right\}
dxdt,
\\
\frac{c}{ 6\log(e+cT) }
\int_0^T \int_{\mathbb{R}^n}
\frac{1}{1+|ct-|x||} 
\left\{
(\partialbar_{c0}u)^2+c^2\sum_{j=1}^n (\partialbar_{cj}u)^2
\right\}
dxdt
\Big]
\\
\le
\int_{\mathbb{R}^n} g^2+c^2|\nabla f|^2 dx
+
2\int_0^T \int_{\mathbb{R}^n} |(\partial_t u) F|dx dt.
\end{multline}
\end{lemma}

\proof 
The estimate for the first term in \eqref{WeightedEnergy-Tangential-n} follows form the standard energy estimate.
The bound for the second term for the case $c=1$ and $n=3$ is given by Lindblad and Rodnianski \cite[p76, Corollary 8.2]{Lindblad-Rodnianski-2005-CMP} 
(see also Lindblad and Rodnianski \cite[p1431, Lemma 6.1]{Lindblad-Rodniansiki-2010-AnnMath}).
We show its generalization following their arguments. 
We put 
\beq
e=(e_0, e_1, \cdots,e_n)
:=
\left(
\frac{1}{2}\{ (\partial_t u)^2+c^2|\nabla u|^2\}, \ -c^2\partial_tu \nabla u
\right).
\eeq
Then we have 
\beq
\label{Div-n}
\partial_te_0+\nabla\cdot (e_1,\cdots,e_n)=(\partial_t u)F.
\eeq
Integrating it on $[0,T]\times \mathbb{R}^n$, we obtain the standard energy estimate 
\beq
\label{ClassicalEnergy-n}
\sup_{0\le t\le T}\int_{\mathbb{R}^n} e_0(t)dx \le 
\int_{\mathbb{R}^n} e_0(0)dx
+
\int_0^T\int_{\mathbb{R}^n} |(\partial_t u)F|dxdt=:X.
\eeq

For $-\infty< q\le T$, we consider the truncated forward light cone 
\beq
\begin{array}{l}
C_0^T(q):=\{(t,x)\ |\ 0\le t\le T,\ |x|= ct-cq\},
\\
K_0^T(q):=\{(t,x)\ |\ 0\le t\le T,\ |x|\le ct-cq\},
\\ 
B_T(q):=\{(T,x)\ |\ |x|\le cT-cq\},
\ \ \ \ 
B_0(q):=\{(0,x)\ |\ |x|\le -cq\}.
\end{array}
\eeq
By the integration of \eqref{Div-n} on $K_0^T(q)$, we have 
\begin{multline}
\label{Y}
Y(q):=\frac{1}{ \sqrt{1+c^2} }
\int_{C_0^T(q)} e\cdot (c, -\frac{x}{r}) d\sigma
\\
=
\int_{B_T(q)} e_0(T) dx -\int_{B_0(q)} e_0(0)dx
-\int\int_{K_0^T(q)} \partial_t uFdxdt
\le 2X,
\end{multline}
where we have used \eqref{ClassicalEnergy-n} for the last inequality.
So that, we have 
\beq
\label{TangentEnergty-Proof-10}
\int_{-\infty}^T \frac{Y(q)}{(1+c|q|)^{1+\kappa} } dq 
\le \int_{-\infty}^T \frac{2X}{(1+c|q|)^{1+\kappa} } dq
\le \frac{4X}{c\kappa}.
\eeq
Since a direct computation shows 
\beq
\label{TangentEnergty-Proof-15}
Y(q)=
\frac{c}{\sqrt{1+c^2}}
\int_{C_0^T(q)} \overline{e}_0 d\sigma,
\eeq
where $\overline{e}_0:= \frac{1}{2}\{(\partialbar_{c0}u)^2+c^2\sum_{j=1}^n (\partialbar_{cj} u)^2\}$, we have 
\beq
\label{TangentEnergty-Proof-20}
\int_{-\infty}^T \frac{Y(q)}{(1+c|q|)^{1+\kappa} } dq 
=
\int_0^T \int_{\mathbb{R}^n} \frac{\overline{e}_0}{(1+|ct-r|)^{1+\kappa} }dxdt, 
\eeq
where we note that $d\sigma dq={\sqrt{1+c^2}}dxdt/{c}$ with 
$q=t-{r}/{c}$.
Combining \eqref{TangentEnergty-Proof-10} and \eqref{TangentEnergty-Proof-20}, we obtain 
\beq
\frac{c\kappa}{2}\int_0^T \int_{\mathbb{R}^n} 
\frac{\overline{e}_0}{(1+|ct-r|)^{1+\kappa} }dxdt 
\le 2X,
\eeq
which is the estimate for the second term in 
\eqref{WeightedEnergy-Tangential-n} as required.

The estimate for the third term in 
\eqref{WeightedEnergy-Tangential-n} follows similarly with slight modification.
From \eqref{Y}, we have 
\beq
\label{TangentEnergty-Proof-100}
\int_{-T}^T \frac{Y(q)}{1+c|q|} dq 
\le \int_{-T}^T \frac{2X}{1+c|q|}dq
\le \frac{4X\log(1+cT)}{c}.
\eeq
By \eqref{TangentEnergty-Proof-15}, we have 
\beq
\label{TangentEnergty-Proof-200}
\int_{-T}^T \frac{Y(q)}{1+c|q|} dq 
=
\int_0^T \int_{t\ge {r}/{c}-T} \frac{\overline{e}_0}{1+|ct-r|}dxdt. 
\eeq
On the other hand, we have 
\beq
\label{TangentEnergty-Proof-300}
\int_0^T \int_{t\le {r}/{c}-T} \frac{\overline{e}_0}{1+|ct-r|}dxdt
\le \int_0^T\frac{dt}{1+cT}
\sup_{0\le t\le T}
\int_{\mathbb{R}^n} \overline{e}_0 dx
\le \frac{2X}{c},
\eeq
where we have used $\overline{e}_0\le 2e_0$ and \eqref{ClassicalEnergy-n} for the last inequality.
Combining \eqref{TangentEnergty-Proof-100}, \eqref{TangentEnergty-Proof-200} and \eqref{TangentEnergty-Proof-300}, we obtain the required estimate.
\ebox

%

\subsection{Wave equations with single speed}

In this subsection, we show that the remark by 
Lindblad and Rodnianski 
for semilinear wave equations in three spatial dimensions 
(see \cite[p52]{Lindblad-Rodnianski-2005-CMP}) 
is also useful for the quasilinear wave equations and the case of two dimensions.

Let $n=2,3$, and let $c>0$, $D\ge1$.
We put $u=(u_1,\cdots,u_D)$, $F=(F_1,\cdots,F_D)$, $f=(f_1,\cdots,f_D)$, $g=(g_1,\cdots,g_D)$, and we consider the Cauchy problem of wave equations with single speed $c$ 
\beq
\label{P-Rn}
\left\{
\begin{array}{l}
(\partial_t^2-c^2\Delta)u_I(t,x)
=
F_I(u',u'')(t,x)
\ \ \mbox{for}\ (t,x)\in [0,\infty)\times \mathbb{R}^n,\ 1\le I\le D 
\\
u(0,\cdot)=f(\cdot),\ \ \partial u(0,\cdot)=g(\cdot),
\end{array}
\right.
\eeq
where we put $\partial_0=\partial_t$ and we denote 
the first derivatives $\{\partial_j u\}_{0\le j\le n}$ by $u'$, and the second derivatives $\{\partial_j\partial_k u\}_{0\le j,k\le n}$ by $u''$.
We assume that $F$ vanishes to the second order when $n=3$, the third order when $n=2$, and has the form 
\begin{equation}
\label{Appendix-F}
F_I(u',u'')=B_I(u')+Q_I(u', u'').
\end{equation}
When $n=3$, $B_I$ and $Q_I$ are given by \eqref{B}, \eqref{Q}, and satisfy the symmetry conditions \eqref{Symmetry}, and the null conditions \eqref{Null-Condition} with $c:=c_1=\cdots=c_D$.
When $n=2$, $B_I$ and $Q_I$ are given by 
\begin{equation}
B_I(u'):= 
\sum_{\substack{1\le J,K,L\le D \\ 0\le j,k,l\le 2} }
B_I^{JKLjkl}\partial_j u_J\partial_k u_K \partial_l u_L
\end{equation}
\begin{equation}
Q_I(u',u''):= 
\sum_{\substack{1\le J,K,L\le D \\ 0\le j,k,l,m\le 2} }
Q_I^{JKLjklm}\partial_j u_J \partial_k u_K \partial_l \partial_m u_L
\end{equation}
and satisfy the symmetry condition
\begin{equation}
Q_I^{JKLjklm}
=
Q_L^{JKIjkml},
\end{equation}
which is required for the energy conservation.
We assume the standard null conditions 
\begin{equation}
\label{Null-Condition-2dim}
\sum_{0\le j,k,l\le 2} B_I^{JKLjkl} \xi_j \xi_k \xi_l
=
\sum_{0\le j,k,l,m\le 2} Q_I^{JKLjklm} \xi_j \xi_k \xi_l \xi_m=0
\end{equation}
for any $1\le I, J, K, L\le D$, and any 
$(\xi_0,\xi_1,\xi_2)\in \mathbb{R}^3$ with   
$\xi_0^2=c^2(\xi_1^2+\xi_2^2)$.
For example, 
$B_I(u')=\sum_{J=1}^D \lambda^J u'\{(\partial_t u_J)^2-c^2|\nabla u_J|^2\}$ 
satisfies the null conditions for $\{\lambda^J\}_{1\le J\le D}\subset \mathbb{R}$, and $Q_I(u',u'')=\partial B_I(u')$ also.

We show an alternative proof of the following theorem.
Lindblad and Rodnianski pointed out the simple proof for the semilinear case with $n=3$ and $c=1$.
We consider the quasilinear case and also the case $n=2$.

%

\begin{theorem}
\label{Theorem-23dim}
Let $n=2$ or $n=3$.
Let $f$ and $g$ be smooth functions.
Then there exist a positive natural number $N$ such that if 
\begin{equation}
\sum_{|\alpha|\le N}\|\anglex^{|\alpha|}\partial_x^\alpha \nabla f\|_{L^2(\Rn)}
+
\sum_{|\alpha|\le N}\|\anglex^{|\alpha|+1}\partial_x^\alpha g\|_{L^2(\Rn)}
=:\varepsilon
\end{equation}
is sufficiently small, then \eqref{P-Rn}
has a unique global solution $u\in C^\infty([0,\infty)\times\Rn)$.
\end{theorem}

We are able to take $N=10$ when $n=3$, and $N=8$ when $n=2$ in the theorem.
The above result for $n=2$ has been shown by 
Godin \cite{Godin-1993-CPDE}, 
Hoshiga \cite{Hoshiga-1995-AdvMathSciAppl},
and Katayama \cite{Katayama-1993-RIMS, Katayama-1995-RIMS}. 
See also Alinhac \cite{Alinhac-2001-InventMath}.
Our proof is based on the weighted energy estimates Lemma \ref{Lemma-WE}, the following Klainerman-Sobolev estimates and the estimates for null conditions.

For any $c>0$, $\Theta_c$ denotes the vector fields 
\begin{equation}
\partial_t,\ \ \partial_j,\ \ 
ct\partial_j+\frac{x_j}{c}\partial_t,\ \ 
x_j\partial_k-x_k\partial_j, 
\ \ 1\le j\neq k\le 3,\ \ \ \ 
t\partial_t+r\partial r
\end{equation}
and $\alpha$ denotes multiple indices.
We note that the $c$ speed Lorentz boosts 
$\{ct\partial_j+{x_j}\partial_t/{c}\}_{j=1}^n$ 
are commutable with $\square_c$, namely, 
$(ct\partial_j+{x_j}\partial_t/{c})\square_c=\square_c(ct\partial_j+{x_j}\partial_t/{c})$.

\begin{lemma}
\label{Lemma-KS-Rn}
(see \cite[p118, Proposition 6.5.1]{Hoermander-1997-Springer}, 
\cite[p43, Theorem 1.3]{Sogge-2008-InterPress})
For any fixed $c>0$, the following estimate holds for any $u$.
\begin{equation}
(1+t+r)^{(n-1)/2} (1+|ct-r|)^{1/2}|u(t,x)|
\lesssim 
\sum_{|\alpha|\le n/2+1} \|\Theta_c^\alpha u(t,\cdot)\|_{L^2(\mathbb{R}^n)}
\end{equation}
\end{lemma}

To estimates the null conditions, we use Lemma \ref{Lemma-NullConditions} with $\Gamma$ replaced by $\Theta_c$ when $n=3$.
When $n=2$, we use the following lemma, which proof is omitted since it is similar to Lemma \ref{Lemma-NullConditions}.
We put $\partialbar_c=(\partialbar_{c0},\partialbar_{c1}, \cdots,\partialbar_{cn}):=(\partial_t+c\partial_r,\nabla-\omega \partial_r)$,
$r=|x|$, $\omega\in \mathbb{S}^{n-1}$, 
$\partial_r=\omega\cdot \nabla$.

\begin{lemma}
\label{Lemma-NullConditions-2dim}
Let $c>0$. 
Let 
\beq
\begin{array}{l}
B(u,v,w)=\sum_{0\le j,k\le 2} B^{jkl}
\partial_ju \partial_kv \partial_lw,
\\
Q(u,v,w)=\sum_{0\le j,k,l\le 2} Q^{jklm}
\partial_ju \partial_k v \partial_l \partial_m w
\end{array}
\eeq 
satisfy the null conditions
\begin{equation}
\sum_{0\le j,k\le 2}
B^{jkl}\xi_j \xi_k \xi_l
=
\sum_{0\le j,k,l,m\le 2}
Q^{jklm}\xi_j \xi_k \xi_l \xi_m=0
\ \ \ \ \mbox{for}\ \ \xi_0^2=c^2(\xi_1^2+\xi_2^2).
\end{equation}
Then the following estimates hold for any $\alpha$ and functions $u$, $v$ and $w$,
where $\beta\le \alpha$ means any component of the multiindices satisfies the inequality.
\begin{multline}
(1) \ \ 
|\Theta_c^\alpha B(u,v, w)|
\lesssim 
\sum_{\substack{
\beta+\gamma+\delta \le\alpha } }
\{
|\partialbar_c \Theta_c^\beta u| |(\Theta_c^\gamma v)'| |(\Theta_c^\delta w)'|
\\
+
|(\Theta_c^\beta u)'| |\partialbar_c \Theta_c^\gamma v| |(\Theta_c^\delta w)'|
+
|(\Theta_c^\beta u)'| |(\Theta_c^\gamma v)'| |\partialbar_c \Theta_c^\delta w|
\}
\end{multline}
\begin{multline}
(2) \ \ 
|\Theta_c^\alpha Q(u,v,w)|
\lesssim 
\sum_{\substack{
\beta+\gamma+\delta\le\alpha } }
\{
|\partialbar_c \Theta_c^\beta u| |(\Theta_c^\gamma v)'| |(\Theta_c^\delta w)''|
\\
+
|(\Theta_c^\beta u)'| |\partialbar_c \Theta_c^\gamma v| |(\Theta_c^\delta w)''|
+
|(\Theta_c^\beta u)'| |(\Theta_c^\gamma v)'| |\partialbar_c (\Theta_c^\delta w)'|
+
\\
+
\frac{1}{\angler} 
|(\Theta_c^\beta u)'||(\Theta_c^\gamma v)'|
(
|(\Theta_c^\delta w)'|
+
|(\Theta_c^\delta w)''|
)
\}
\end{multline}
\end{lemma}

%

The proof of Theorem \ref{Theorem-23dim} is similar to that of Theorem \ref{Theorem}.
It suffices to prove the following proposition.
In its proof, we implicitly use $\varepsilon^2\le C\varepsilon$ since $\varepsilon$ is sufficiently small.

\begin{proposition}
\label{Prop-23dim}
Let $n=2$ or $n=3$.
Let $M_0$ be a positive number which satisfies
$M_0\ge 6$ when $n=3$, 
$M_0\ge 4$ when $n=2$.
We put 
\beq
\varepsilon:=
\sum_{|\alpha|\le M_0+4}
\|\anglex^{|\alpha|} \partial^\alpha \nabla f\|_{L^2(\mathbb{R}^n)} 
+
\sum_{|\alpha|\le M_0+4}
\|\anglex^{|\alpha|+1}\partial^{\alpha} g\|_{L^2(\mathbb{R}^n)}.
\eeq
Let $T>0$ and $A_0>0$. 
We put $S_T:=[0,T)\times \Rn$.
Let $u\in C^\infty([0,T)\times \mathbb{R}^n)$ be the local solution of \eqref{P-Rn}.
We assume 
\begin{equation}
\label{Ineq-Assumption-23dim}
\sum_{\substack{|\alpha|\le M_0 } }
\sup_{(t,x)\in S_T
}
\angletpr^{(n-1)/2}\anglectmr^{1/2}
|\Theta_c^\alpha u'(t,x)|
\le A_0 \varepsilon.
\end{equation}
If $\varepsilon$ is sufficiently small, 
then there exist constants $C_0>0$, which is independent of $A_0$, and $C>0$, which is dependent on $A_0$, such that the following estimates hold.
\beq
(1)\ \ \ \ 
\sum_{|\alpha|\le M_0+4} 
\|\Theta_c^\alpha u'(t,\cdot)\|_{L^2(\Rn)}
\le C\varepsilon (1+t)^{C\varepsilon }
\ \ \mbox{for}\ \ 0\le t< T
\eeq
\begin{multline}
(2)\ \ \ \ 
\sum_{|\alpha|\le M_0+3}
(\log(e+t))^{-1/2}
\|\anglex^{-1/2}\Theta_c^\alpha u'\|_{L^2(S_t)}
\\
+
\sum_{\substack{|\alpha|\le M_0+3} }
(\log(e+t))^{-1/2}
\|\langle cs-r\rangle^{-1/2} 
\partialbar_c \Theta_c^\alpha u\|_{L^2(S_t)}
\\
\le
C\varepsilon (1+t)^{C\varepsilon}
\ \ \mbox{for}\ \ 0\le t< T
\end{multline}
\begin{multline}
(3)\ \ \ \ 
\sum_{\substack{|\alpha|\le M_0+2} }
\sup_{(s,x)\in S_t}
\langle s+r \rangle^{(n-1)/2}\langle cs-r \rangle^{1/2}
|\Theta_c^\alpha u'(s,x)|
\\
\le
C\varepsilon (1+t)^{C\varepsilon}
\ \ \mbox{for}\ \ 0\le t< T
\end{multline}
\beq
(4)\ \ \ \ 
\sum_{|\alpha|\le M_0+2 }
\|\Theta_c^\alpha u'\|_{L^\infty((0,T),L^2(\Rn))}
\le C_0\varepsilon + C\varepsilon^{(6-n)/2}
\eeq
\begin{multline}
(5)\ \ \ \ 
\sum_{\substack{|\alpha|\le M_0} }
\sup_{(t,x)\in S_T
}
\angletpr^{(n-1)/2}\anglectmr^{1/2}
|\Theta_c^\alpha u'(t,x)|
\le
C_0\varepsilon + C\varepsilon^{(6-n)/2}
\end{multline}
\end{proposition}

\proof
Put $M=M_0+4$.
First, we remark that under the assumption 
\eqref{Ineq-Assumption-23dim}, we have  
\begin{equation}
\label{Ineq-Assumption-Es-23dim}
\begin{array}{ll}
\displaystyle
(1+t+r)
\sum_{|\alpha|\le M/2+1}|\Theta_c^\alpha u'(t,x)|\le C \varepsilon
& 
\mbox{when}\ \ n=3 \\
\displaystyle
(1+t+r)
\sum_{|\alpha|\le M/3+1}|\Theta_c^\alpha u'(t,x)|^2\le C \varepsilon^2
& 
\mbox{when}\ \ n=2
\end{array}
\end{equation}
for some constant $C>0$ since $M/(5-n)+1\le M_0$.

(1) 
The proof is essentially same to that of Proposition \ref{Prop} by the use of \eqref{Ineq-Assumption-23dim}.
For $1\le I,L\le D$ and $0\le l,m\le n$, we put 
\beq
\gamma_I^{Llm}
:=
\sum_{ \substack{1\le J,K\le D \\ 0\le j,k\le n} }
Q_I^{JKLjklm}\partial_ju_J\partial_k u_K.
\eeq
For any $\alpha$ with $|\alpha|\le M$, we use \eqref{Eqn-Divergence} and its integration to have 
\begin{multline}
\partial_t\int_{ \mathbb{R}^n } e_0(\Theta_c^\alpha u) dx
\le 
C
\sum_{1\le I\le D}
\|\square_{ \gamma_I } \Theta_c^\alpha u_I\|_{L^2(\mathbb{R}^n)}
\|(\Theta_c^\alpha u)'\|_{ L^2(\mathbb{R}^n) }
\\
+
C\sum_{ \substack{1\le I,L\le D \\ 0\le l,m\le n} }
\|\partial_{t,x}\gamma_I^{Llm}\|_{ L^\infty(\mathbb{R}^n) }
\|(\Theta_c^\alpha u)'\|_{ L^2(\mathbb{R}^n) }^2. 
\end{multline}
Similarly to \eqref{Ineq-Commutator}, we have 
\begin{multline}
\sum_{\substack{1\le I\le D \\  |\alpha|\le M} }
\|\square_{ \gamma_I } \Theta_c^\alpha u_I\|_{L^2(\mathbb{R}^n)}
\lesssim
\sum_{\substack{1\le I\le D \\  |\alpha|\le M} }
\| \Theta_c^\alpha \square_{ \gamma_I }u_I\|_{L^2(\mathbb{R}^n)}
\\
+
\sum_{\substack{|\alpha|+|\beta|\le M \\ |\beta|\le M-1} }
\sum_{\substack{1\le I,L\le D \\ 0\le l,m\le n} }
\|(\Theta_c^\alpha \gamma_I^{Llm}) \Theta_c^\beta u''\|_{L^2(\mathbb{R}^n)}
\lesssim 
\frac{\varepsilon^{4-n}}{1+t} \sum_{|\alpha|\le M}\|\Theta_c^\alpha u'\|_{L^2(\mathbb{R}^n)},
\end{multline}
where we have used \eqref{Ineq-Assumption-Es-23dim} for the last inequality.
Since 
$\sum_{|\alpha|\le M}\|\Theta_c^\alpha u'\|_{L^2(\mathbb{R}^n)}$ 
is equivalent to  
$\sum_{|\alpha|\le M} \{\int_{\mathbb{R}^n} e_0(\Theta_c^\alpha u)dx\}^{1/2}$ 
for small $\varepsilon$, 
we obtain 
\beq
\sum_{\substack{|\alpha|\le M} }
\partial_t 
\left\{\int_{ \mathbb{R}^n } e_0(\Theta_c^\alpha u) dx\right\}^{1/2}
\le 
\frac{C\varepsilon^{4-n}}{1+t} 
\sum_{|\alpha|\le M}
\left\{\int_{ \mathbb{R}^n } e_0(\Theta_c^\alpha u) dx\right\}^{1/2},
\eeq
which leads to the required inequality by the Gronwall inequality. 

%

(2)
By Lemma \ref{Lemma-EnergyKSSTangential}, we have 
\begin{multline}
\label{Ineq-Linear-M}
\sum_{|\alpha|\le M-1}
(\log(e+t))^{-1/2}\|\anglex^{-1/2} \Theta_c^\alpha u'\|_{L^2(S_t)}
\\
+
\sum_{\substack{|\alpha|\le M-1} }
(\log(e+t))^{-1/2}
\|\langle cs-r\rangle^{-1/2} \partialbar_c \Theta_c^\alpha u\|
_{L^2(S_t)} 
\\
\le
C_0\varepsilon
+
C_0\sum_{\substack{|\alpha|\le M-1} }
\|\Theta_c^\alpha \square_{c} u\|_{L^1((0,t),L^2(\mathbb{R}^n))}.
\end{multline}
The last term is bounded by 
\begin{multline}
\sum_{\substack{|\alpha|\le M-1} }
\|\Theta_c^\alpha \square_{c} u\|_{L^1((0,t),L^2(\mathbb{R}^n))}
\le
C\int_0^t \frac{\varepsilon^{4-n}}{1+s} \sum_{|\alpha|\le M}\|\Theta_c^\alpha u'(s, \cdot)\|_{L^2(\mathbb{R}^n)} ds
\\
\le
C\varepsilon^{4-n} (1+t)^{C\varepsilon},
\end{multline}
where we have used \eqref{Ineq-Assumption-Es-23dim}, (1) and  
\begin{equation}
\label{Ineq-Nonlinear-Base}
\sum_{\substack{ |\alpha|\le M-1 } }
|\Theta_c^\alpha \square_{c} u|
\lesssim 
\left(
\sum_{|\beta|\le M/(5-n)}|\Theta_c^\beta u'|
\right)^{4-n}
\sum_{|\alpha|\le M}|\Theta_c^\alpha u'|.
\end{equation}

%

(3) The estimate follows from 
Lemma \ref{Lemma-KS-Rn} and (1).

%

(4) By the standard energy estimates, we have 
\begin{multline}
\sum_{|\alpha|\le M_0+2}
\|\Theta_c^\alpha u_I'(t,\cdot)\|_{L^2(\Rn)}^2
\le 
C_0\sum_{|\alpha|\le M_0+2}
\|\Theta_c^\alpha u_I'(0,\cdot)\|_{L^2(\Rn)}^2
\\
+
C_0\sum_{|\alpha|\le M_0+2}
\int_0^t\int_{\Rn} 
|\partial_t \Theta_c^\alpha u_I \square_{c} \Theta_c^\alpha u_I|dxds
=: A_1+A_2.
\end{multline}
We have $A_1\le (C_0\varepsilon)^2$ for some $C_0>0$ which is independent of $A_0$.
By Lemma \ref{Lemma-NullConditions} with $\Gamma$ replaced by $\Theta_c$ and Lemma \ref{Lemma-NullConditions-2dim}, 
$A_2$ is bounded as $A_2\lesssim A_3+A_4$, 
where
\begin{multline}
A_3:=
\int_0^t\int_\Rn
\left(
\sum_{|\alpha|\le M_0+2} |\Theta_c^\alpha u'|
\right)^{4-n}
\sum_{|\alpha|\le M_0+3} |\partialbar_{c} \Theta_c^\alpha u|
\sum_{|\alpha|\le M_0+3} |\Theta_c^\alpha u'|
dxds
\end{multline}
\begin{multline}
A_4:=
\int_0^t\int_\Rn
\left(
\sum_{|\alpha|\le M_0+2} |\Theta_c^\alpha u'|
\right)^{4-n}
\sum_{|\alpha|\le M_0+3} |\Theta_c^\alpha u'|
\sum_{|\alpha|\le M_0+3} |\Theta_c^\alpha u'|
\frac{dx}{\angler}ds.
\end{multline}
We use (3) to have 
\begin{multline}
A_3\le 
C\varepsilon^{4-n}
\sum_{|\alpha|\le M_0+3}
\|\langle cs-r\rangle^{-1/2}\langle s \rangle^{-\delta} \partialbar_{c} \Theta_c^\alpha u\|_{L^2(S_t)}
\\
\cdot
\sum_{|\alpha|\le M_0+3}
\|\angler^{-1/2} \langle s\rangle^{-\delta} \Theta_c^\alpha u'\|_{L^2(S_t)}
\le C\varepsilon^{6-n},
\end{multline}
where $\delta>0$ is a sufficiently small number and we have used 
(2) to obtain the last inequality.
Similarly, we have 
\begin{multline}
A_4\le 
C\varepsilon^{4-n}
\sum_{|\alpha|\le M_0+3}
\|\angler^{-1/2}\langle s\rangle^{-\delta} \Theta_c^\alpha u'\|_{L^2(S_t)}
\\
\cdot
\sum_{|\alpha|\le M_0+3}
\|\angler^{-1/2}\langle s\rangle^{-\delta} \Theta_c^\alpha u'\|_{L^2(S_t)}
\le C\varepsilon^{6-n}.
\end{multline}
Combining these estimates, we obtain the required result.

(5) The estimate follows from Lemma \ref{Lemma-KS-Rn} and (4).
\ebox

%

\end{document}